\documentclass[aos,seceqn,nameyear,dvips]{arximspdf}
\usepackage[mathscr]{euscript}
\usepackage{multirow}
\usepackage{graphicx}

\doi{10.1214/10-AOS845}
\volume{39}
\issue{1}
\pubyear{2011}
\firstpage{362}
\lastpage{388}

\makeatletter

\newtheorem{theorem}{Theorem}[section]
\newtheorem{lem}{Lemma}[section]

\newcommand{\iint}{\int\!\!\int}

\makeatother

\begin{document}
\begin{frontmatter}

\title{Functional single index models for longitudinal~data}
\runtitle{Functional SIMs}

\begin{aug}
\author[A]{\fnms{Ci-Ren} \snm{Jiang}\corref{}\ead[label=e1]{crjiang@stat.berkeley.edu}} and
\author[B]{\fnms{Jane-Ling} \snm{Wang}\thanksref{m2}\ead[label=e2]{wang@wald.ucdavis.edu}}
\runauthor{C.-R. Jiang and J.-L. Wang}
\affiliation{University of California, Berkeley
and University of California, Davis}
\address[A]{Department of Statistics\\
University of California, Berkeley\\
367 Evans Hall\\
Berkeley, California 94720-3860 \\
USA\\
\printead{e1}}
\address[B]{Department of Statistics\\
University of California, Davis\\
4118 Mathematical Sciences Building\\
One Shields Avenue\\
Davis, California 95616\\
USA\\
\printead{e2}}
\end{aug}

\thankstext{m2}{Supported in part by NSF Grant DMS-09-06813.}

% HISTORY:
\received{\smonth{12} \syear{2009}}
\revised{\smonth{7} \syear{2010}}

% ABSTRACT
%
\begin{abstract}
A new single-index model that reflects the time-dynamic effects of
the single index is proposed for longitudinal and functional
response data, possibly measured with errors, for both
longitudinal and time-invariant covariates. With appropriate
initial estimates of the parametric index, the proposed estimator
is shown to be $\sqrt{n}$-consistent and asymptotically normally
distributed. We also address the nonparametric estimation of
regression functions and provide estimates with optimal
convergence rates. One advantage of the new approach is that the
same bandwidth is used to estimate both the nonparametric mean
function and the parameter in the index. The finite-sample
performance for the proposed procedure is studied numerically.
\end{abstract}

% KEYWORDS
%

\end{frontmatter}

%s1 ###
\section{Introduction}

For univariate response variables $Y$ with multivariate covariate
$Z\in\mathbb{R}^p$, the single-index model
%
%e1.1 ###
\begin{equation}
\label{s1} \mathbb{E}(Y| Z) = m(\beta_0^TZ)
\end{equation}
is an attractive dimension-reduction method to model the effect of
multivariate covariates nonparametrically. Since $m(\cdot)$, known as
the \textit{link function}, is an unknown smooth function, the scale of
$\beta_0^TZ$ may be determined arbitrarily. For identifiability
reasons, $\beta_0$ is often assumed to be a unit vector with
nonnegative first coordinate. The primary parameter of interest is
the coefficient $\beta_0$ in the index $\beta_0^T Z$ since
$\beta_0$ makes explicit the relationship between the response variable
$Y$ and the covariate $Z$. There are several different approaches to
estimate $\beta_0$ in (\ref{s1}), such as the projection pursuit
regression [\citet{FrieS81}, \citet{Hall89}], average derivatives
[\citet{HardS89}, \citet{Ichi93}] and partial least-squares
[\citet{NaikT00}] methods. Typically, the link function needs to be
undersmoothed in order to estimate $\beta_0$ at the $\sqrt{n}$-rate.
\citet{HardHI93} showed that a $\sqrt{n}$-consistent estimator of
$\beta_0$ can be achieved without undersmoothing the link function,
that is, the same bandwidth can be used to estimate both the parameter
$\beta_0$ and the nonparametric link function $m(\cdot)$. However,
their approach relies on a grid search to obtain the estimate for
$\beta_0$ and is time consuming when the dimension $p$ is high. To
overcome this drawback, and inspired by the sliced inverse regression
method [\citet{Li91}, \citet{XiaTLZ02}] proposed a new method,
called ``conditional minimum average variance estimation''
(MAVE). Unlike most previous methods, MAVE does not need to
undersmooth the nonparametric link function estimator to attain the
$\sqrt{n}$-rate consistency for the parametric index estimate. Also,
it does not require strong assumptions on the distribution of the
covariates. Theoretical results for this approach to single-index
models are available in \citet{Xia06} and some extensions have
been studied in \citet{Xia07} and \citet{KongX07}, among
others. However, none
of these works addresses longitudinal data, which is the focus of
this paper.

Our goal is to extend MAVE to the following single-index models
for functional/longitudinal response data:
%
%e1.2 ###
\begin{equation}\label{s2}
\mathbb{E}(Y(t)| Z(t)) = \mu(t, \beta_0^TZ(t)),
\end{equation}
where $Y(t)$, $t\in\mathcal{T}$, is a stochastic process on a
compact time interval $\mathcal{T}$, $Z$ contains $p$
covariates, some or all of which may be stochastic functions over
the time interval $\mathcal{T}$, and, to be identifiable, $\beta_0$ is
a unit vector with
nonnegative first coordinate. More specifically,
%
%e1.3 ###
\begin{equation}\label{s3}
Y(t)= \mu(t, \beta_0^TZ(t))+ \epsilon(t, Z(t)),
\end{equation}
where $\mu$ is an unknown bivariate link function and
$\epsilon(t, Z(t))$ is a random function with mean 0 that reflects
the within-subject correlations of measurements and possibly
measurement errors at different time points. Thus, there are two
distinctive features in the functional single-index model
(\ref{s2}), as compared to the traditional single-index model
$(\ref{s1})$ considered in \citet{Xia06} and \citet{XiaTLZ02}.
First, the functional single-index model accommodates longitudinal
response and longitudinal covariates, as well as vector covariates.
Second, the effects of the single index and, consequently,
covariates $Z$, may change over the time dynamic through a
bivariate link function and this seems more realistic for
longitudinal responses.

Recently, \citet{BaiFZ09} combined penalized splines and quadratic
inference functions to estimate the index coefficient and unknown
link function in a single-index model for longitudinal data.
However, the link function in their model is univariate and thus
does not reflect the dynamic effects of the single index. Moreover,
their approach is restricted to generalized linear models, where the
variance function of the response is a known function of the mean
function. In contrast, the link function in our model is an unknown
function of time and the index, reflecting the dynamic feature of
the effect of the single index, and the structure of the variance
function is not restricted in our approach.

The rest of this paper is organized as follows. Section \ref{sec2} extends
the original MAVE method to longitudinal data. Asymptotic theory
for the proposed estimators is described in Section \ref{sec3}, with proofs
in the \hyperref[appB]{Appendix}. Practical implementations of the new approaches
and simulation studies are presented in Section \ref{sec4}. In Section \ref{sec5},
we apply our method to two AIDS data sets: one with time-invariant
covariate and the other also involving longitudinal covariates.
Section \ref{sec6} contains our conclusions.

%s2 ###
\section{Methodology}\label{sec2}

We begin with the setting of model (\ref{s1}) for univariate
response $Y$ and multivariate covariate $Z \in\mathbb{R}^p$.
Let $\sigma_\beta(\beta^TZ)$ be the conditional variance of $Y$
given $\beta^TZ$. The true direction $\beta_0$ in (\ref{s1}) is the
solution of $\beta$ that minimizes
$\mathbb{E}\{\sigma_\beta(\beta^TZ)\}=\mathbb{E}\{Y-\mathbb{E}(Y|\beta
^TZ)\}^2$.

For a random sample, $\{(Y_i,Z_i), i=1,\ldots,n \}$, of $(Y,Z)$,
$\mathbb{E}(Y|\beta^TZ)$ can be approximated locally at $\beta^T
Z_j$ by a linear expansion, that is, $\mathbb{E}(Y|\beta^TZ)\approx a_j+
b_j^T\beta^T(Z-Z_j)$. Empirically, $\sigma_\beta(\beta^TZ)$ can be
approximated at $\beta^T Z_j$ by
$\sum_{i=1}^n[Y_i-\{a_j+b_j^T\beta^T(Z_i-Z_j)\}]^2w_{ij}$, where
$w_{ij}\geq0$ are weights with $\sum_{i=1}^n w_{ij}=1$, for example, $
w_{ij}=K_h\{\beta^T(Z_i-Z_j)\}/\sum_{i=1}^nK_h\{\beta^T(Z_i-Z_j)\}$,
where $K_h(\cdot)=h^{-d}K(\cdot/h)$ and $d$ is the dimension
of $K(\cdot)$. Therefore, we can estimate $\beta_0$ by solving the
minimization problem
%
%e2.1 ###
\begin{equation}\label{rM1}
\min_{\beta,\mathbf{a},\mathbf{b}}\Biggl(\sum_{j=1}^n\sum_{i=1}^n
[Y_i-\{a_j+b_j^T\beta^T(Z_i-Z_j)\}]^2w_{ij} \Biggr),
\end{equation}
where $\mathbf{a}=(a_1,\ldots,a_n)$ and $\mathbf{b}=(b_1,\ldots,b_n)$. Given
$\beta$, (\ref{rM1}) is a local linear smoother of the data
$\{Y_i,\beta^T(Z_i-Z_j)\}$, while, given $\mathbf{a}$ and $\mathbf{b}$,
(\ref{rM1}) is just a weighted least-squares problem for $\beta$.
Consequently, the minimization in (\ref{rM1}) can be viewed as a
combination of nonparametric function estimation and parametric
direction estimation.
Furthermore, the weights can be updated iteratively via the
relation
$\tilde{w}_{ij}=K_h\{\hat{\beta}^T(Z_i-Z_j)\}/\sum_{i=1}^nK_h\{\hat
{\beta}^T(Z_i-Z_j)\}$,
using the current estimate $\hat{\beta}$, then updating the
estimate of $\beta_0$ by minimizing (\ref{rM1}) with $w_{ij}$
replaced by $\tilde{w}_{ij}$. This could be repeated until
$\hat{\beta}$ converges and is called refined MAVE (rMAVE) in
\citet{XiaTLZ02}.

%s2.1 ###
\subsection{Estimation}
Hereafter, the response will be longitudinal data, which typically
consists of random fluctuations or measurement errors. Let
$Y_{ij}=Y_i(T_{ij})$ be the $j$th observation for the $i$th subject,
made at a random time $T_{ij} \in\mathcal{T}$, where $\mathcal{T}$
is an interval. Along with the responses, we have information on $p$
covariates, some of which may be longitudinal covariates. Since a
univariate covariate can be considered a special case of
a longitudinal covariate with constant value, we will adopt the
notation for longitudinal covariates and define $Z_{ij}=Z_i(T_{ij})
\in\mathbb{R}^p$, $i=1,\ldots,n$ and $j=1,\ldots,N_i$, as the
$p$-dimensional covariate for the $i$th subject evaluated at time
$T_{ij}$. The functional single-index model (\ref{s3}) applied to
the observed longitudinal and covariates data leads to
\[ %\begin{equation}\label{m1}
Y_{ij} = \mu(T_{ij},\beta_0^TZ_{ij})+ \epsilon(T_{ij},Z_{ij}).
\]
For simplicity, we only consider bounded covariates $Z$ when
deriving theoretical properties, even though our simulation study
shows that the method could work well for unbounded covariates. The
bounded assumption is commonly adopted in the literature, for example,
in \citet{HardHI93} and \citet{HardS89}.
Here, we assume that the measurement times $T_{ij}$ are a random
sample of size $N_i$, assumed to be i.i.d. and independent
of all other random variables.

The two main steps in our approach are to estimate the direction
$\beta_0$ and the mean function $\mu$. In particular, we show how to
estimate the parametric index $\beta_0$ by adapting rMAVE for
longitudinal data. The asymptotic distribution of $\hat{\beta}$ is
studied in Section \ref{sec3} for both longitudinal and time-invariant
covariates. The mean function can then be estimated through a
two-dimensional scatter plot smoother of $Y_{ij}$ on
$(T_{ij},\hat{\beta}^TZ_{ij})$ when $\hat{\beta}$ is available.

To estimate the parametric index efficiently, we extend rMAVE to
longitudinal data. For simplicity, and to avoid the curse of
dimensionality, we only consider a single index in our model.
Therefore, $\beta$ is a vector instead of a matrix. As in MAVE,
for any given $(T_{j\ell},Z_{j\ell})$,
$\mathbb{E}(Y_{ik}|T_{ik},\beta^TZ_{ik})$ can be approximated by a
linear expansion at $(T_{j\ell}, \beta^TZ_{j\ell})$, that is,
$\mathbb{E}(Y_{ik}|T_{ik},\beta^TZ_{ik})\approx a_{j\ell} +
b_{j\ell}(T_{ik}-T_{j\ell})+d_{j\ell}\beta^T(Z_{ik}-Z_{j\ell})$.
Similarly, the conditional covariance,
$\sigma_\beta(T_{ik},\beta^TZ_{ik})=\mathbb{E}\{Y_{ik}-\mathbb
{E}(Y_{ik}|T_{ik},\beta^TZ_{ik})\}^2$,
can be approximated by
$\sum_{i=1}^n\sum_{k=1}^{N_i}[Y_{ik}-\{a_{j\ell}+
b_{j\ell}(T_{ik}-T_{j\ell})+d_{j\ell}\beta^T(Z_{ik}-Z_{j\ell})\}
]^2w_{ikj\ell}$,
where
%
%e2.2 ###
\begin{eqnarray}\label{wt}
w_{ikj\ell} &=& \frac{K(({T_{ik}-T_{j\ell}})/{h_t},({\beta
^T(Z_{ik}-Z_{j\ell})})/{h_z})}
{\sum_{i=1}^n\sum_{k=1}^{N_i}K(({T_{ik}-T_{j\ell}})/{h_t},({\beta
^T(Z_{ik}-Z_{j\ell})})/{h_z})},\hspace*{-25pt}\nonumber\\[-8pt]\\[-8pt]
\sum_{i=1}^n\sum_{k=1}^{N_i} w_{ikj\ell} &=& 1.\nonumber
\end{eqnarray}
Here, $K(\cdot)$ is a two-dimensional kernel function of order $(0,2)$
defined in Appendix \ref{appC} with compact
support that is also a symmetric density function with finite
moments of all orders and bounded derivatives; $h_t$ and $h_z$
are the respective bandwidths for smoothing along the time ($t$) and
single-index covariate ($\beta^T z$) direction. We can then estimate
$\beta_0$ by solving the minimization problem
%
%e2.3 ###
\begin{eqnarray}\label{rM2}
&&\min_{\beta,\mathbf{a},\mathbf{b},\mathbf{d}}\Biggl(\sum_{j=1}^n\sum_{\ell=1}^{N_j}
\sum_{i=1}^n\sum_{k=1}^{N_i}
[Y_{ik}-\{a_{j\ell}+b_{j\ell}(T_{ik}-T_{j\ell})\nonumber\\[-8pt]\\[-8pt]
&&\qquad\hspace*{107.6pt}{}+d_{j\ell}\beta
^T(Z_{ik}-Z_{j\ell})\}]^2w_{ikj\ell}
\Biggr).\nonumber
\end{eqnarray}
Suppose that we have a current estimator $\hat{\beta}$ of $\beta_0$
and current refined weights $\tilde{w}_{ikj\ell}$.
The estimate for $\beta_0$ will be updated by minimizing equation
(\ref{rM2}) with $w_{ikj\ell}$ replaced by $\tilde{w}_{ikj\ell}$.
This procedure will be repeated until $\hat{\beta}$ converges. The
final estimate, $\hat{\beta}$, can then be used to estimate the mean
function $\mu$ via a two-dimensional smoother that has the
same bandwidth as the weights in (\ref{rM2}), that is,
%
%e2.4 ###
\begin{eqnarray}
\label{eqMu}
&\hspace*{-10pt}\displaystyle \hat{\mu}(t,\hat{\beta}^Tz) = \hat{b}_0\qquad \mbox{where, for
} \mathbf{b}=(b_0,b_1,b_2),&\nonumber\\
&\hspace*{-10pt}\displaystyle \hat{\mathbf{b}} = \mathop{\arg\min}_{\mathbf{b}} \sum_{i=1}^{n}\sum_{j=1}^{N_i}
K\biggl\{\frac{t-T_{ij}}{h_t},\frac{\hat{\beta}^T(z-Z_{ij})}{h_z}\biggr\}&
\\[-3pt]
&\hspace*{-10pt}\displaystyle \hspace*{52.4pt}\hspace*{63.4pt}{} \times
\{Y_{ij}-b_0-b_1(T_{ij}-t)-b_2\hat{\beta}^T(Z_{ij}-z)\}^2.\hspace*{-52.4pt}&\nonumber
\end{eqnarray}

%s2.2 ###
\subsection{Algorithm}\label{sec22}

Let $h_t$ and $h_z$ be the bandwidths for $T$ and $\beta^TZ$,
respectively, and let $\hat{\sigma}_{\beta}^2$ denote the quantity
to be minimized in (\ref{rM2}), which is within the parentheses.
Define $K_h(t,z)=K(t/h_t,z/h_z)/(h_th_z)$.

\begin{enumerate}
\item Start with an initial value of $\beta$, say
$\hat{\beta}_{(0)}$.

\item Use the current estimate $\hat{\beta}_{(m)}$ and weighted
least-squares method to obtain
$(\hat{\mathbf{a}},\hat{\mathbf{b}},\hat{\mathbf{d}})=\arg\min_{\mathbf{a},\mathbf{b},\mathbf{d}}
\hat{\sigma}_{\beta_{(m)}}^2$, where
\begin{eqnarray*}
w_{ikj\ell}&=&K_h \bigl\{
(t_{ik}-t_{j\ell}),\hat{\beta}_{(m)}^T(z_{ik}-z_{j\ell}) \bigr\}\\
&&{}\Big/
\sum_{i=1}^{n}\sum_{k=1}^{N_i}K_h \bigl\{
(t_{ik}-t_{j\ell}),\hat{\beta}_{(m)}^T(z_{ik}-z_{j\ell}) \bigr\}.
\end{eqnarray*}

\item Use the estimates $(\hat{\mathbf{a}},\hat{\mathbf{b}},\hat{\mathbf{d}})$
from step 2 to obtain the updated estimate
$\hat{\beta}_{(m+1)}=\arg\min_{\beta}\hat{\sigma}_\beta^2$.

\item Repeat steps 2 and 3 until
$\|\hat{\beta}_{(m+1)}-\hat{\beta}_{(m)}\|<\varepsilon$, where
$\varepsilon$ is some given tolerance value.

\item The final estimate of $\beta$ from step 4 is then used to
reach the final estimate of the mean function defined in
(\ref{eqMu}).
\end{enumerate}

%s2.3 ###
\subsection{Bandwidth selection}

Instead of selecting the bandwidths by the leave-one-curve-out
cross-validation method suggested in \citet{RiceS91}, we choose
the bandwidths for the mean function estimator via an $m$-fold
cross-validation procedure to reduce the computational cost. Below,
we describe the $m$-fold cross-validation method for the bandwidth
selection for $\mu(t,\beta^Tz)$. Supposing that subjects are
randomly divided into $m$ groups, ($S_1, S_2,\ldots,S_m$), the
$m$-fold cross-validation bandwidth is
%
%e2.5 ###
\begin{equation} \label{eqkfm} h_\mu=\mathop{\arg\min}_{h}\sum_{\ell=1}^{m}
\sum_{i\in S_\ell}\sum_{j=1}^{N_i}
\bigl\{Y_{ij}-\hat{\mu}^{(-S_\ell)}(T_{ij},\hat{\beta}^TZ_{ij})\bigr\}^2,
\end{equation}
where $\hat{\mu}^{(-S_\ell)}(T_{ij},\hat{\beta}^TZ_{ij})$ is the
estimated mean function at $(T_{ij},\hat{\beta}^TZ_{ij})$,
excluding subjects in $S_\ell$.

%s3 ###
\section{Asymptotic results}\label{sec3}

We assume that $(T_{ij},Z_{ij},Y_{ij})$ have the same distribution as
$(T,Z,Y)$ with joint probability density function $g_3(t,z,y)$ and that
the observational times
$T_{ij}$ are i.i.d. with probability density function $g(t)$, but
dependency is
allowed among observations from the same subject. Let
$\tilde{z}=\beta^Tz$, $\tilde{z}^0=\beta_0^Tz$ and
let $f_2(t,\tilde{z})$ and $f_3(t,\tilde{z},y)$ be the joint densities
of $(T, \tilde{Z})$ and $(T, \tilde{Z}, Y)$, respectively.
%The following derivations are under Assumptions A.1-A.6, listed in
%Appendix A, and
The kernel function is assumed to be symmetric. For simplicity, we
also assume that $\int u^2K(u,v)=\int v^2 K(u,v) = \int u^2v^2K(u,v)
=1$ as, without loss of generality, any symmetric density kernel
function can be applied
after proper normalization. Since we are interested in the
asymptotic distribution of $\hat{\beta}$, similar to the assumption
in \citet{HardHI93}, we assume that the initial value
$\hat{\beta}_{(0)}$ is in a $\sqrt{n}$-neighbor of $\beta_0$. This
assumption is for technical convenience; in the simulations, an
arbitrary initial value was used and it performed well. To be prudent, one
may want to try different random initial $\hat{\beta}_{(0)}$ and choose the
final estimate as the one that leads to the smallest value in the
minimization problem of (\ref{rM2}). In the data analysis, we chose ten
different initial values for $\hat{\beta}_{(0)}$ and they all
converged to the same estimate $\hat{\beta}$.

From the iterative algorithm in Section \ref{sec22}, the updated
$\hat{\beta}$ from minimizing (\ref{rM2}) after one iteration will
become
%
%e3.2 ###
%e3.1 ###
\begin{eqnarray*}
\hat{\beta} & = & \beta_0 + \{D_n^{\beta}\}^{-1}\Upsilon+
o_p(n^{-1/2})\qquad\mbox{where } \\[-25pt]
\end{eqnarray*}
\begin{eqnarray}\label{beq1}\qquad
D_n^\beta & = &
\{n\mathbb{E}N\}^{-2}\sum_{i=1}^n\sum_{k=1}^{N_i}\sum_{j=1}^n\sum_{\ell=1}^{N_j}
\frac{d_\beta^2(T_{ik},Z_{ik})}{\hat{f_2}(T_{ik},\tilde
{Z}_{ik})}\nonumber\\[-8pt]\\[-8pt]
&&\hspace*{106.5pt}{} \times
K_h\{(T_{j\ell}-T_{ik}),
(\tilde{Z}_{j\ell}-\tilde{Z}_{ik})\}\nonumber\\
&&\hspace*{106.5pt}{} \times(Z_{j\ell}-Z_{ik})(Z_{j\ell
}-Z_{ik})^T,\nonumber
\\
\label{Upeq0}
\Upsilon & = & (n\mathbb{E}N )^{-2}
\sum_{i=1}^{n}\sum_{k=1}^{N_i} \sum_{j=1}^n\sum_{\ell=1}^{N_j}
\frac{d_{\beta}(T_{j\ell},Z_{j\ell})}{\hat{f}_2(T_{j\ell},\tilde
{Z}_{j\ell})}\nonumber\\
&&\hspace*{106.5pt}{} \times K_h\{(T_{ik}-T_{j\ell}),(\tilde{Z}_{ik}-\tilde{Z}_{j\ell})\}
(Z_{ik}-Z_{j\ell})\nonumber\\
&&\hspace*{106.5pt}{} \times\{Y_{ik}-a_\beta(T_{j\ell},Z_{j\ell})\\
&&\hspace*{122.3pt}{}  -
b_\beta(T_{j\ell},Z_{j\ell})(T_{ik}-T_{j\ell})\nonumber\\
&&\hspace*{122.3pt}{}  -
d_\beta(T_{j\ell},Z_{j\ell})(\tilde{Z}^0_{ik}-\tilde{Z}^0_{j\ell})\},\nonumber
\end{eqnarray}
%
%Through some tedious calculations (Proofs of \ref{Deq1} and
%approximations of $\{D_n^\beta\}^{-1}$ and $\Upsilon$:
$\hat{f}_2(t,z)$ is the estimate of $f_2(t,z)$ and $a_\beta$'s, $b_\beta
$'s and $d_\beta$'s are the coefficients of the linear approximation,
as defined in Section \ref{sec2}. By means of some tedious calculations
[sketches of proofs of
(\ref{Deq1}) and (\ref{Upeq1}) are in Appendix \ref{appB} with assumptions
A.1--A.6 listed in Appendix \ref{appA}], we can obtain the following
approximations of $\{D_n^\beta\}^{-1}$ and $\Upsilon$:
%
%e3.3 ###
\begin{eqnarray}\label{Deq1}\qquad
\{D_n^\beta\}^{-1} &=& \frac{\beta_0\beta_0^T}{\tau}
-\frac{h_z^2}{2\tau}(\tilde{G}^+
\tilde{F}^T\beta_0^T+\beta_0\tilde{F} \tilde{G}^+) +
\frac{1}{2}\tilde{G}^+ + O_p(h+\delta_\beta),
\\
%
%
%e3.4 ###
\label{Upeq1}
\Upsilon &=& \tilde{G}(\beta-\beta_0)-(n\mathbb{E}N
)^{-1} \sum_{i=1}^n\sum_{k=1}^{N_i} \biggl\{
\nu_{\beta_0}(T_{ik},Z_{ik})\,\frac{\partial\mu}{\partial
\tilde{z}^0} \biggr\}\epsilon_{ik}\nonumber\\[-8pt]\\[-8pt]
&&{} +o_p(n^{-1/2}),\nonumber
\end{eqnarray}
where $\tilde{G}=\mathbb{E}\{(\partial\mu/\partial
\tilde{z}^0)^2G(Z)\}/2$, $\tilde{G}^+ =
B_0(B_0^T\tilde{G}B_0)^{-1}B_0^T$ is the Moore--Penrose inverse of
$\tilde{G}$ with $(\beta_0, B_0)$ an orthogonal matrix,
$\tau=\mathbb{E}\{(\partial
\mu/\partial\tilde{z}^0)^2\}h_t^2$,
$\tilde{F}=\mathbb{E}\{(\partial\mu/\partial
\tilde{z}^0)^2F_\beta(T,Z)\}$, $F_\beta(t,z) =
\frac{\partial}{\partial
\tilde{z}}\{f_2(t,\tilde{z})\nu_\beta^T(t,z)\}/f_2(t,\tilde{z})$,
$G(z)=\mathbb{E}\{(Z_{ik}-z)(Z_{ik}-z)^T\}$,
$\nu_{\beta_0}(t,z)=\mathbb{E}(Z|T=t,\beta_0^TZ=\beta_0^Tz)-z$
and $\delta_\beta=|\hat{\beta}-\beta_0|$.

After plugging (\ref{Deq1}) and (\ref{Upeq1}) into (\ref{beq1}),
we obtain
\begin{eqnarray*}
\hat{\beta} & = & \beta_0 + \biggl\{ \frac{\beta_0\beta_0^T}{\tau}
-\frac{h_z^2}{2\tau}(\tilde{G}^+
\tilde{F}^T\beta_0^T+\beta_0\tilde{F} \tilde{G}^+) +
\frac{1}{2}\tilde{G}^+ + O_p(h+\delta_\beta) \biggr\} \\
&&{} \times\Biggl[ \tilde{G}(\beta-\beta_0)-(n\mathbb{E}N
)^{-1} \sum_{i=1}^n\sum_{k=1}^{N_i} \biggl\{
\nu_{\beta_0}(T_{ik},Z_{ik})\,\frac{\partial\mu}{\partial
\tilde{z}^0}
\biggr\}\epsilon_{ik} +o_p(n^{-1/2}) \Biggr] \\
&&{} + o_p(n^{-1/2})\\
& = & \beta_0 (1+c_n) +
\frac{1}{2}(I-\beta_0\beta_0^T)(\beta-\beta_0)-
\frac{\tilde{G}^+}{2n\mathbb{E}N} \sum_{i=1}^n\sum_{k=1}^{N_i}
\biggl\{ \nu_{\beta_0}(T_{ik},Z_{ik})\,\frac{\partial\mu}{\partial
\tilde{z}^0} \biggr\}\epsilon_{ik}\\
&&{} +o_p(n^{-1/2}),
\end{eqnarray*}
where $ c_n= h_z^2\tilde{F} \tilde{G}^+
[(n\mathbb{E}N )^{-1} \sum_{i=1}^n\sum_{k=1}^{N_i}
\{ \nu_{\beta_0}(T_{ik},Z_{ik})(\partial\mu/\partial
\tilde{z}^0) \}\epsilon_{ik} -
\tilde{G}(\beta-\beta_0)]/(2\tau)$.

Since $|\beta|=1$, $\hat{\beta}$ needs to be standardized. From
the above calculation, $|\hat{\beta}|=1+c_n+o_p(n^{-1/2})$ so
\begin{eqnarray*}
\frac{\hat{\beta}}{|\hat{\beta}|}&=&\beta_0 +
\frac{1}{2}(I-\beta_0\beta_0^T)(\beta-\beta_0)\\
&&{}-
\frac{\tilde{G}^+}{2n\mathbb{E}N} \sum_{i=1}^n\sum_{k=1}^{N_i}
\biggl\{ \nu_{\beta_0}(T_{ik},Z_{ik})\,\frac{\partial\mu}{\partial
\tilde{z}^0} \biggr\}\epsilon_{ik} +o_p(n^{-1/2}).
\end{eqnarray*}

Therefore, in the $(m+1)$th iteration,
\begin{eqnarray*}
\hat{\beta}_{(m+1)} & = & \beta_0 +
\frac{1}{2}(I-\beta_0\beta_0^T)\bigl(\hat{\beta}_{(m)}-\beta_0\bigr)\\
&&{}-
\frac{\tilde{G}^+}{2n\mathbb{E}N} \sum_{i=1}^n\sum_{k=1}^{N_i}
\biggl\{ \nu_{\beta_0}(T_{ik},Z_{ik})\,\frac{\partial\mu}{\partial
\tilde{z}^0} \biggr\}\epsilon_{ik} +o_p(n^{-1/2}) \\
& = & \beta_0 +
\frac{1}{2^m}(I-\beta_0\beta_0^T)\bigl(\hat{\beta}_{(1)}-\beta_0\bigr)\\
&&{} -
\Biggl(\sum_{j=1}^m\frac{1}{2^j}\Biggr)(n\mathbb{E}N)^{-1}\tilde{G}^+
\sum_{i=1}^n\sum_{k=1}^{N_i} \biggl\{
\nu_{\beta_0}(T_{ik},Z_{ik})\,\frac{\partial\mu}{\partial
\tilde{z}^0} \biggr\}\epsilon_{ik} \\
&&{} + o_p(n^{-1/2}).
\end{eqnarray*}
Consequently, as the iteration $m\rightarrow\infty$, Lemma
\ref{lem1} in Appendix \ref{appD} implies the following theorem.
\begin{theorem}\label{theo31}
Under assumptions \textup{A.1--A.6},
\[
\sqrt{n}(\hat{\beta}-\beta_0)\stackrel{\mathscr{D}}{\rightarrow}
N_p(0,\Sigma),
\]
where $\Sigma= \tilde{G}^+\Sigma^*\tilde{G}^+$ and
$\Sigma^*= \frac{\mathbb{E}(N)-1}{\mathbb{E}N} \mathbb{E}[\{
\frac{\partial\mu}{\partial\tilde{z}^0}\, \nu_{\beta_0}(T,Z)\epsilon
\}\{\frac{\partial\mu}{\partial\tilde{z}^0}\,
\nu_{\beta_0}(T',Z')\times\break\epsilon'\}^T] +
\frac{1}{\mathbb{E}N}\mathbb{E}[\{ \frac{\partial
\mu}{\partial\tilde{z}^0} \,\nu_{\beta_0}(T,Z)\epsilon
\}\{\frac{\partial\mu}{\partial\tilde{z}^0}\,
\nu_{\beta_0}(T,Z)\epsilon\}^T]$.
\end{theorem}

In practice, the covariance of $\hat{\beta}$ in Theorem \ref{theo31} is
unknown and needs to be estimated to make inference on $\beta$. The
idea is to replace the unknown values with consistent estimates.
First, $\mathbb{E}(N)$ can be estimated by $\bar{N}=\sum_{i=1}^n
N_i/n$ and $\tilde{G}$ can be estimated by
\[
\hat{\hspace*{-1.5pt}\tilde{G}}=\sum_{i=1}^n\sum_{k=1}^{N_i}
\biggl\{\frac{\partial}{\partial
\tilde{z}}\,\hat{\mu}(T_{ik},\hat{\beta}^TZ_{ik})\biggr\}^2\hat
{G}(Z_{ik})/(2n\bar{N}),
\]
where $\hat{G}(Z_{ik})=\frac{1}{n}\sum_{j=1}^n
\frac{1}{N_i}\sum_{\ell=1}^{N_i}(Z_{j\ell}-Z_{ik})(Z_{j\ell}-Z_{ik})^T$.
To estimate $\Sigma^*$, we estimate $\nu_{\beta_0}(T,Z)$ at all
$(T_{ik},Z_{ik})$ by (\ref{nueq}), estimate $\epsilon$ at
$(T_{ik},Z_{ik})$ by the residual,
$Y_{ik}-\hat{\mu}(T_{ik},\hat{\beta}^TZ_{ik})$, and average the
product terms in $\Sigma^*$. Therefore,
\[
\hat{\Sigma}^* =
\frac{\bar{N}-1}{\bar{N}}\Biggl\{\frac{1}{N^*}\sum_{i=1}^n\sum_{1\leq
j \neq k \leq N_i}H_{ik} H_{ij}^T\Biggr\} +
\frac{1}{\bar{N}}\Biggl\{\frac{1}{n\bar{N}}\sum_{i=1}^n\sum
_{k=1}^{N_i}H_{ik}H_{ik}^T\Biggr\},
\]
where $H_{ik}=\frac{\partial}{\partial
\tilde{z}}\,\hat{\mu}(T_{ik},\hat{\beta}^TZ_{ik})\hat{\nu}_{\beta
_0}(T_{ik},Z_{ik})\hat{\epsilon}_{ik}$
and $N^*=\sum_{i=1}^nN_i^2-N_i$. To estimate
$\nu_{\beta_0}(T_{ik},Z_{ik})=
\mathbb{E}(Z|T=T_{ik},\beta_0^TZ=\beta_0^TZ_{ik})-Z_{ik}$, we can,
for simplicity, apply a weighted average estimator on the
observations in the neighborhood of $(T_{ik}, \beta_0^TZ_{ik})$,
which leads to
%
%e3.5 ###
\begin{equation}\label{nueq}\quad
\hat{\nu}_{\beta_0}(T_{ik},Z_{ik}) = \sum_{j,\ell}
\frac{K_h\{(T_{j\ell}-T_{ik}),\hat{\beta}^T(Z_{j\ell}-Z_{ik})\}}{\sum
_{j,\ell}K_h\{(T_{j\ell}-T_{ik}),\hat{\beta}^T(Z_{j\ell}-Z_{ik})\}}
Z_{j\ell} - Z_{ik}.
\end{equation}
More sophisticated procedures might be considered to estimate the
above unknown values.

Before showing the asymptotic property of the local linear
smoother, $\hat{\mu}(t,\break\hat{\beta}^Tz(t))$, we first need the
asymptotic property of the local linear smoother,
$\hat{\mu}(t,u(t))$, where $u$ is a univariate longitudinal
covariate. This is provided in Theorem \ref{thm02} below, with the proof in
Appendix \ref{appC}.
\begin{theorem}\label{thm02}
Under assumptions \textup{A.1--A.6}, $h_z/h_t\rightarrow\rho$ and $
n\mathbb{E}(N)h_t^{6}\rightarrow\tau^2$ for some $0< \rho, \tau< \infty
$ and
\[
\sqrt{n\bar{N}h_th_z}[\hat{\mu}(t,u)-\mu(t,u)]
\stackrel{\mathscr{D}}{\rightarrow}N( \eta(t, u),\Sigma_{\mu} (t,u)),
\]
where $\eta(t, u)=\frac{\tau\sqrt{\rho}}{2}
\{
\frac{\partial^2 \mu}{\partial t^2} + \frac{\partial^2 \mu}{\partial
u^2}\rho^2
\}$, $\Sigma_\mu(t,u) =
[\operatorname{var}(Y|t,u)\|K_2\|^2]/f_2(t,u)$,\break $\|K_2\|^2=\int K_2^2$
and $f_2(t,u)$ is the joint density of $(T,U)$.
\end{theorem}

It is interesting that the asymptotic bias term in Theorem \ref{thm02}
depends on the ratio $h_z/h_t$.
This is due to the assumption that $nE(N)h_t^{6}\rightarrow\tau^2$ and
assumption A.1 in Appendix \ref{appA}, which requires $h_t$ and $h_z$ to have
the same rate. After some Taylor expansions, these two assumptions on
the bandwidths lead to the asymptotic bias term in Theorem \ref{thm02}. The
assumption that the two bandwidths, $h_z$ and $h_t$, have the same rate is
natural since the mean function $\mu(t,z)$ has the same order of
smoothness along both the $t$ and $z$ coordinates.

Since $\hat{\beta}$ is a $\sqrt{n}$-consistent estimator of
$\beta_0$ and by Theorem \ref{thm02}, the asymptotic properties of
the local linear smoothers for the mean can be obtained.
\begin{theorem}
Under assumptions \textup{A.1--A.6}, $h_z/h_t\rightarrow\rho$ and $
n\mathbb{E}(N)h_t^{6}\rightarrow\tau^2$ for some $0< \rho, \tau< \infty$, and
\[
\sqrt{n\bar{N}h_th_z}\{\hat{\mu}(t,\hat\beta^Tz)-\mu(t,\beta_0^Tz)\}
\stackrel{\mathscr{D}}{\rightarrow} N(\triangle_\mu,\Sigma_\mu),
\]
where $\triangle_\mu=\frac{\tau\sqrt{\rho}}{2}
\{
\frac{\partial^2 \mu}{\partial t^2} + \frac{\partial^2 \mu}{\partial
\tilde{z}^2}\rho^2
\}$, $\Sigma_\mu =
[\operatorname{var}(Y|t,\beta_0^Tz)\|K_2\|^2]/f_2(t,\beta_0^Tz)$ and
$\|K_2\|^2=\int K_2^2$.
\end{theorem}

%s4 ###
\section{Simulation study}\label{sec4}

Two simulation schemes are considered in this paper. One considers
the case with time-invariant covariates; the other considers
longitudinal covariates. In both simulation studies,
$\hat{\beta}_{(0)}=(1/\sqrt{p})\mathbf{1}_p$, where $\mathbf{1}_p$ is a
$p$-dimensional vector with entries all equal to
1, was used as the initial value of $\beta$, the number of runs was
100 and the number of subject for each run was $n=100$.

%s4.1 ###
\subsection{Simulation \textup{I}: Time-invariant covariate}
The covariate for each subject $(Z_1, Z_2)$ is generated from
$Z_1 \sim \operatorname{Bernoulli}$ (with probability of success 0.5) and $Z^T_2
\sim N_5(0,\Psi)$, where $\Psi=0.5\times I_5+ 0.5\times\mathbf{1}_5\mathbf{1}_5^T$.
We choose $\beta_0^T =(2, 1, 0, 3, 0, -1)/\sqrt{15}$.
%and $\tilde{Z}=(Z_1,Z_2^T)\beta_0$.
Given a subject with covariate $\tilde{z}=(z_1, z_2)^T \beta_0$, the
stochastic process $Y^*$ is generated from a Gaussian process on
$[0, 1]$ with mean function $
\mu(t,\tilde{z})=\sin(\tilde{z})\sin(t\pi)+\{1-\sin(\tilde{z})\}\cos
(t\pi)$
and covariance function $\Gamma(s, t)=(1/4)\phi_1(t)\phi_1(s) +
(1/16)\phi_2(t)\phi_2(s)$, where $\phi_1(t)=-\sin(\pi t)\sqrt{2}$
and $\phi_2(t)=\cos(\pi t)\sqrt{2}$ are the eigenfunctions of
$\Gamma$, with corresponding eigenvalues $1/4$ and $1/16$,
respectively. The measurement errors are assumed to be normally
distributed with mean $0$ and variance $0.01$. Note that the
variance of measurement error is not very small compared to the
scales of the mean function and the eigenvalues.

For the measurement schedule, we use a ``jittered'' design with an
equally spaced grid $\{c_0,\ldots,c_{50}\}$ on $[0,1]$ ($c_0=0$
and $c_{50}=1$) and then jitter each point $c_i$ by
$s_i=c_i+\epsilon_i$, where $\epsilon_i$ are i.i.d. with
$N(0,0.0001)$, $s_i=0$ if $s_i<0$ and $s_i=1$ if $s_i>1$. This
resulted in a jittered schedule that is no longer equally spaced;
from there, a random sample of size $N_i$ is selected from
$\{s_1,\ldots,s_{49}\}$ without replacement to serve as the $N_i$
measurement schedule for the $i$th subject, where $N_i$ is itself
sampled from a discrete uniform distribution $\{2,\ldots,10\}$.

%
%t1 ###
\begin{table}
\caption{Performances of estimators for 3- and 10-fold
CV. (Measures of differences between $\hat{\beta}$ and $\beta$:
\mbox{$\| \cdot \|$} measures the\vspace*{1pt} difference in the Euclidean norm and
$\cos^{-1}$ in terms of the angle between $\hat{\beta}$ and
$\beta$.) The IMSE\vspace*{1pt} is the integrated mean squared error, defined as
$\iint[\hat{\mu}(t,u)-\mu(t,u)]^2\,dt\,du$. Note that the IMSE of
$\hat{\mu}$ in simulation \textup{II} is much larger than that in
simulation \textup{I}, due to different scales of $t$ and $\beta^Tz$, and
different mean function} \label{tb_imseI}
\begin{tabular*}{\tablewidth}{@{\extracolsep{\fill}}lcccc@{}}
\hline
\textbf{Simulation} & \textbf{CV} & $\bolds{\|\beta-\hat{\beta}\|}$ & $\bolds{\cos^{-1}(\beta^T\hat
{\beta})}$ &\textbf{IMSE}$\bolds{(\mu)}$ \\
\hline
\phantom{I}I & \phantom{0}3 & 0.2121 (0.0867)& 12.1900 (5.0161) & 0.0312 \\
& 10 & 0.2121 (0.0793) & 12.1860 (4.5833)& 0.0237 \\
[2pt]
II & \phantom{0}3 & 0.2575 (0.1923)& \phantom{0}14.8944 (11.3004) & 0.2723
\\
&10 & 0.2501 (0.1914) & \phantom{0}14.4675 (11.2666) & 0.2257 \\
\hline
\end{tabular*}
\end{table}

%
%f1 ###
\begin{figure}[b]

\includegraphics{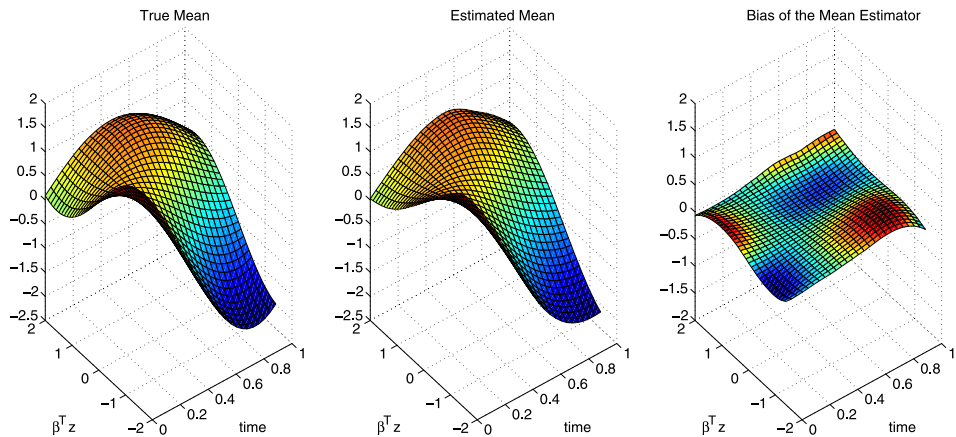}

\caption{True mean function, averaged estimated mean
function and bias in simulation study \textup{I}.}
%The average
%estimate for the mean function and its bias show the estimating
%procedures work well when 10-fold CV is employed to select the
%bandwidth. The results when 3-fold CV is employed are similar and
%to save the space is not provided here.\fi}
\label{S4M10}
\end{figure}

We experimented with several $m$-fold cross-validation (CV) methods and
found the 10-fold method to be satisfactory. Table \ref{tb_imseI} reports
the results for $m=3$ and~$10$. The results for the parametric
estimate $\hat{\beta}$ are comparable for $3$- and $10$-fold CVs
with the 10-fold CV being somewhat better. In terms of estimating the
mean function, the $10$-fold CV performs better. Figure \ref{S4M10} also
suggests good performance of the $10$-fold CV method in terms of
bias. The plot for the $3$-fold CV is similar, but is not provided.

%s4.2 ###
\subsection{Simulation \textup{II}: Longitudinal covariates}

This simulation scheme is inspired by the CD4$+$ cell counts data
from the Multicenter AIDS Cohort Study or MACS [\citet{KaslODP87}],
which is analyzed in Section \ref{sec5}. There are five covariates in this
AIDS data: age at seroconversion and four longitudinal covariates
[packs of cigarettes, recreational drug use (1: yes, 0: no),
number of sexual partners and mental illness scores (CESD), larger
values indicate increased depressive symptoms]. In the simulation,
the covariate values were based on the five covariates from 100
randomly selected subjects. The mean function, coefficient of the
index, two eigenfunctions and two eigenvalues were also chosen to
mimic the corresponding values of the real data (see Section \ref{sec5}).
Therefore, we choose the mean function
\[
\mu(t,\beta_0^Tz)=6+\frac{\beta_0^Tz}{5}+\frac{1}{1+\exp(t)}+
\frac{\exp\{-t(\beta_0^Tz+3)\}}{1+\exp\{-t(\beta_0^Tz+3)\}},
\]
where the index coefficient $\beta_0^T=(0.1043, 0.5213, 0.8341,
-0.1043, -0.1043)$ and $t\in[-3,5.5]$. The two eigenfunctions are
$\phi_1(t)=\cos\{(t+3)\pi/8.5\}/\sqrt{4.25}$, and
$\phi_2(t)=-\sin\{(t+3)\pi/8.5\}/\sqrt{4.25}$, with respective
eigenvalues $\lambda_1=2$ and $\lambda_2=0.5$. For each subject,
the two principal component scores are generated from
$N(0,\lambda_1)$ and $N(0,\lambda_2)$. Also, normally distributed
measurement errors with mean zero and variance $0.1$ are added.

%
%f2 ###
\begin{figure}

\includegraphics{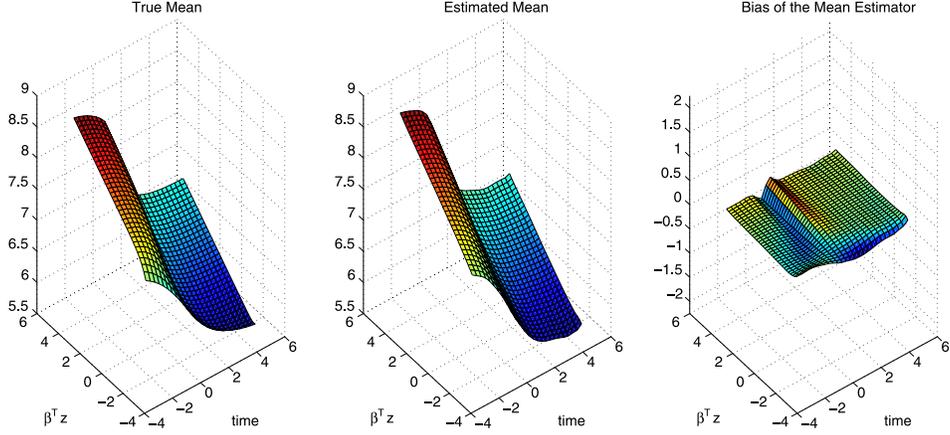}

\caption{True mean function, averaged estimated mean
function and bias in simulation study \textup{II}.} %\iffalse The average
%estimate for the mean function and its bias show the estimating
%procedures work well when 10-fold CV is employed to select the
%bandwidth. The results when 3-fold CV is employed are similar and
%to save the space is not provided here.\fi}}
\label{S5M10}
\end{figure}

Consistent with the results in simulation I, where the covariates
are time invariant, the results in Table \ref{tb_imseI} for
estimating $\beta$ and mean function $\mu(t,\beta^Tz(t))$ are also
comparable for the 3-fold and 10-fold CVs, with 10-fold CV slightly
better. Again, we only provide the plot of estimates based on the
10-fold CV. Other than the boundary, Figure \ref{S5M10} suggests good
performance of the 10-fold CV method in terms of bias. The boundary
effect appears to be due to the sparsity of the data and is more
prominent than in the previous simulation setting of time-invariant
covariates. The observed $\beta^Tz$ are very sparse near the
boundaries.

%s5 ###
\section{Application}\label{sec5}

We illustrate the methodology via CD4$+$ cell counts data from
the Multicenter AIDS Cohort Study or MACS [\citet{KaslODP87}]. HIV
destroys CD4 cells, which play a vital role in the immune system.
The CD4 cell count is thus a good marker for disease progress. The
number of CD4 cells might also be related to some subject-specific
factors such as smoking, age, etc. In the first example, we apply
our approach to the CD4 data analyzed in \citet{WuC00}, where the
covariates are time invariant. The second example is the CD4 data
analyzed in \citet{ZegeD94}, where longitudinal covariate
variables are available.

%s5.1 ###
\subsection{Example \textup{I}: Time-invariant covariates}

This data set involves
1817 measurements of CD4 percentages, which are cell counts divided
by the total number of lymphocytes, observed for 283 homosexual
men who became HIV positive between 1984 and 1991. The measurements
were scheduled at each half-yearly visit; however, the actual
measurement times may vary and some subjects missed some of their
scheduled visits. The resulting measurement times $t_{ij}$ per
subject are irregular and sparse. Three time-independent covariate
variables are considered in our analysis: smoking status
(1: yes, 0: no), age at HIV infection and pre-HIV infection CD4
percentage. To make the scales of different covariates compatible,
we standardize age and pre-HIV infection CD4 percentage. Similarly to
the simulation study, we use 3-fold and 10-fold CV to choose the
bandwidths for the nonparametric procedures.

To avoid\vspace*{2pt} being trapped in a local minimum, we choose ten different
random initial values for $\hat{\beta}_{(0)}=
(\hat{\beta}_{1(0)},\hat{\beta}_{2(0)},\hat{\beta}_{3(0)})^T$, as
follows. First, we pick five different values (0.1, 0.3, 0.5, 0.7
and 0.9) for $\hat{\beta}_{1(0)}$ and generate
$\hat{\beta}_{2(0)}$ from $U(0,\sqrt{1-\hat{\beta}_{1(0)}^2})$,
then we set
$\hat{\beta}_{3(0)}=\sqrt{1-\hat{\beta}_{1(0)}^2-\hat{\beta}_{2(0)}^2}$
to ensure that $\|\hat{\beta}_{(0)}\|=1$. The signs of
$\hat{\beta}_{2(0)}$ and $\hat{\beta}_{3(0)}$ are initially
randomly assigned and then flipped to make up for the ten initial
$\hat{\beta}_{(0)}$. These ten initial values all lead to the same
$\hat{\beta}$.

%
%t2 ###
\begin{table}[b]
\caption{Estimated parametric index $\hat{\beta}$ and
asymptotic covariance of $\hat{\beta}$ for example \textup{I}\break [here,
$h_{\mu}=(h_t,h_z)$ is the bandwidth for estimating $\beta$ and
$\mu(t,\beta^Tz)$]}\label{cd4bV}
\begin{tabular*}{\tablewidth}{@{\extracolsep{\fill}}lcc@{\hspace*{-2pt}}}
\hline
CV & 3 & 10 \\
[4pt]
$h_\mu$ & $(2.14, 2.80)$ & $(1.70, 5.00)$ \\[4pt]
$\hat{\beta}^T$ & $(0.0727, -0.1074, 0.9916)$ & $(0.0877, -0.1076,
0.9903)$ \\
[4pt]
$\operatorname{Var}(\hat{\beta})\approx\frac{\hat{\Sigma}}{\sqrt{n}}$ &
$
\pmatrix{0.4213 & 0.0141 & -0.0796 \cr
0.0141 & 0.0887 & -0.0161 \cr
-0.0796 & -0.0161 & 0.0602}
$ & $
\pmatrix{
0.4137 & 0.0103 & -0.1072 \cr
0.0103 & 0.0898 & -0.0128 \cr
-0.1072 & -0.0128 & 0.0932}$ \\
\hline
\end{tabular*}
\end{table}

%
%f3 ###
\begin{figure}

\includegraphics{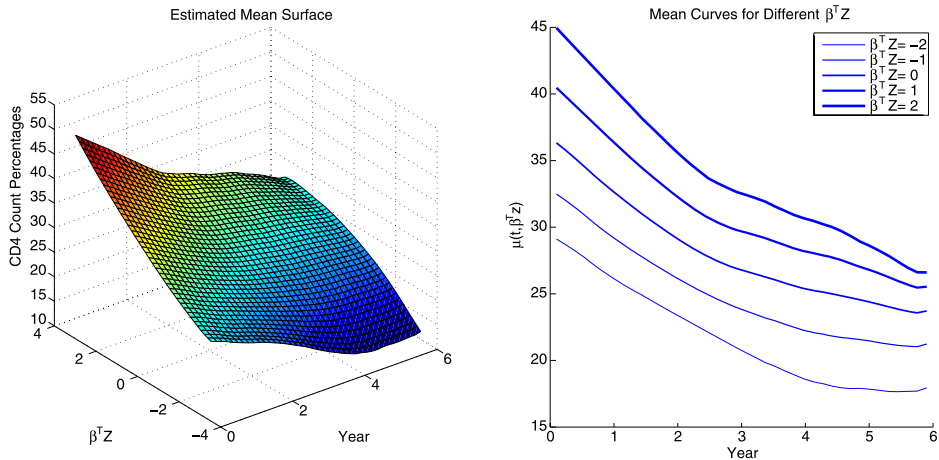}

\caption{Estimated mean function of AIDS data for example
\textup{I}.} %\iffalse The average CD4 percentages deplete quickly at the
%beginning of HIV infection and the time when the rate of depletion
%slows down varies with $\hat{\beta}^Tz$. The rate of depletion
%slows down earlier when the subject has higher $\hat{\beta}^Tz$.\fi}}
\label{cd4mean}
\end{figure}

Several statistical models have been applied to this data set, such
as varying coefficient models in \citet{WuC00}. In their analysis,
only the effect of pre-infection CD4 percentage was found to be
significant and positive, but none of the covariate effects seem
time-dependent [see Figures 1 and 2 in \citet{WuC00}]. This
result is consistent with our findings in Table \ref{cd4bV} and
Figure \ref{cd4mean}. We find that people who smoke, who are young
when they get the HIV infection and who have higher pre-HIV
infection CD4 percentages tend to have higher post-HIV infection CD4
percentages on average. However, only pre-HIV infection CD4
percentage is significant.
In the right panel of Figure \ref{cd4mean}, we observe that, in
general, the CD4 percentages deplete rather quickly at the
beginning of HIV infection and the rates of depletion during the
first 2.5 years are generally higher than in later years.
However, the time when the rate of depletion slows down varies
with the levels of $\hat{\beta}^Tz$. More specifically, when
$\hat{\beta}^Tz$ is larger, the rate of depletion slows down
earlier.

%s5.2 ###
\subsection{Example \textup{II}: Time-invariant and longitudinal covariates}
%%%%%%%%% Longitudinal Covariate AIDS data %%%%%%%%%%%%%%%%%%%%%%%

In this data set, 2376 CD4 observations on 369 subjects were made
and the times of observation ranged from 3 years before to 6 years
after seroconversion. Five covariates considered in this analysis
are age, packs of cigarettes, recreational drug use (1: yes, 0:
no), number of sexual partners and mental illness scores (CESD)
(larger values indicate increased depressive symptoms). Except for
age, the other four covariates are longitudinal. As in example I,
we applied $3$- and $10$-fold CVs to choose the bandwidths in
nonparametric procedures and adopted the same strategy to select 10
initial values for $\hat{\beta}_0$. It turned out that all ten
random initial $\hat{\beta}_{(0)}$'s lead to the same $\hat{\beta}$.

Previous analysis for this data includes the semiparametric models
in \citet{ZegeD94}, where
%Some statistical models have been applied to this data set, such as
%semiparametric models in
age, smoking, recreational drug
use and increased numbers of sexual partners are associated with
higher CD4 cell numbers, while increased depressive symptoms are
associated with decreased CD4 levels, but the effects of age and
recreation drug use are not significant.

%
%t3 ###
\begin{table}
\tablewidth=315pt
\caption{Estimated parametric index $\hat{\beta}$ and
asymptotic covariance of $\hat{\beta}$ for example \textup{II}\break [here,
$h_{\mu}=(h_t,h_z)$ is the bandwidth for estimating $\beta$ and
$\mu(t,\beta^Tz)$]}\label{cd4bV2}
\begin{tabular*}{\tablewidth}{@{\extracolsep{\fill}}lc@{\hspace*{-2pt}}}
\hline
\multicolumn{2}{c}{\textbf{3-fold CV}}\\
%CV & 3 & 10 \\ \hline
\hline
$h_\mu$ & $(1.25, 3.00)$ \\[2pt]
$\hat{\beta}^T$ & $(0.0141, 0.5700, 0.8211, -0.0159, -0.0216)$ \\
[3pt]
$\operatorname{Var}(\hat{\beta})\approx\frac{\hat{\Sigma}}{\sqrt{n}}$ &
$
\pmatrix{0.0035 & 0.0045 & 0.0210 & -0.0010 & -0.0003 \cr
0.0045 & 0.0956 & 0.2733 & -0.0049 & -0.0038 \cr
0.0210 & 0.2733 & 2.4311 & 0.0029 & -0.0214 \cr
-0.0010 & -0.0049 & 0.0029 & 0.0069 & -0.0009 \cr
-0.0003 & -0.0038 & -0.0214 & -0.0009 & 0.0021}
$ \\
\hline
\multicolumn{2}{c}{\textbf{10-fold CV}}\\
\hline
$h_\mu$ & $(1.00, 4.00)$ \\[2pt]
$\hat{\beta}^T$ & $(0.0128, 0.5530, 0.8326, -0.0193, -0.0225)$ \\[3pt]
$\operatorname{Var}(\hat{\beta})\approx\frac{\hat{\Sigma}}{\sqrt{n}}$ & $
\pmatrix{0.0037 & 0.0070 & 0.0284 & -0.0011 & -0.0005 \cr
0.0070 & 0.1287 & 0.3744 & -0.0065 & -0.0061 \cr
0.0284 & 0.3744 & 2.7206 & -0.0018 & -0.0302 \cr
-0.0011 & -0.0065 & -0.0018 & 0.0070 & -0.0008 \cr
-0.0005 & -0.0061 & -0.0302 & -0.0008 & 0.0023}$
\\
\hline
\end{tabular*}
\end{table}

%
%f4 ###
\begin{figure}[b]

\includegraphics{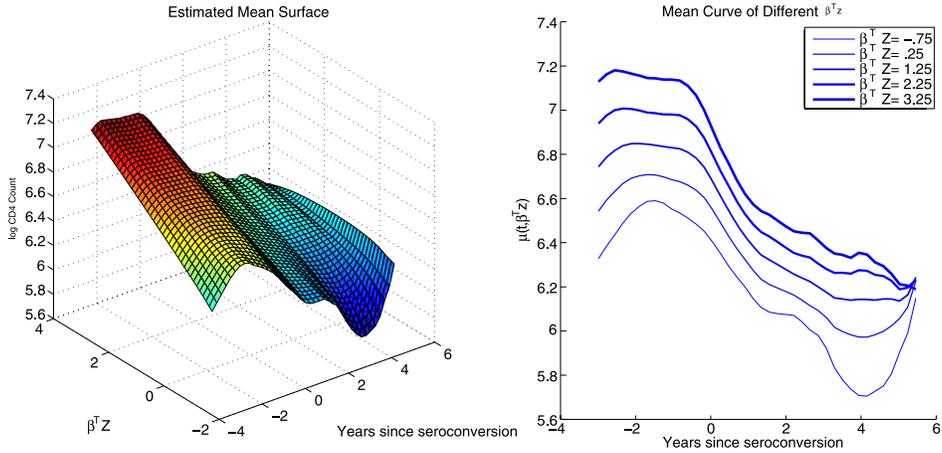}

\caption{Estimated mean function of AIDS data for example \textup{II}.}
%at the time of seroconversion and the time when the rate of
%depletion slows down slightly varies with $\hat{\beta}^Tz(t)$. \fi
\label{cd4mean2}
\end{figure}

In our analysis, among these five covariates, the effect of packs
of cigarettes smoked per day is the most significant. Moreover,
our analysis in Table \ref{cd4bV2} suggests that an increasing number of
sexual partners is negatively associated with CD4 counts, which
seems more reasonable than the previous result.
%The later could be because the immune system of an AIDS person
%with more sexual partners is weaker on the average.
From Table~\ref{cd4bV2} and Figure \ref{cd4mean2}, we also observe
higher mean CD4 cell numbers when subjects are older, smoke more,
use recreational drugs and have lower CESD. The right panel of
Figure \ref{cd4mean2} suggests a big decline of CD4 cell counts
half a year before the seroconversion and one year after
seroconversion. After one year of seroconversion, the decline
slows down slightly. The trough at the end might be due to the
boundary effect.

%s6 ###
\section{Discussion and conclusions}\label{sec6}

The proposed estimate of the single index for longitudinal response
data works well in our simulations and data application. It has the
advantage that the same bandwidth is used to estimate the
nonparametric mean function and the single index without the need to
undersmooth the mean function in order to achieve the $\sqrt
{n}$-convergence rate, as is often the case in semiparametric regression
models with independent response data. This leads to a unified
approach to selecting the bandwidth. Additional computational savings
are accomplished through the $m$-fold cross-validation methods. The
simulation results reported in Section \ref{sec4} suggest that the
performance of the procedure is not very sensitive to the choice of
$m$ and the initial value $\hat{\beta}_{(0)}$.

We have derived asymptotic distributions for both the parametric
($\beta$) and nonparametric ($\mu$) components of the model and
illustrate its usefulness for statistical inference via an AIDS
data set. While it is possible to extend the approach to multiple
indices, such an extension would be computationally intensive and
subject to the curse of dimensionality.

An additive model
\[
\mathbb{E}(Y(t))=\mu(t, \beta_0^TZ(t)) = \mu_t(t) + \mu_z(\beta_0^TZ(t))
\]
can be viewed as a special case of model (\ref{s3}), and if the model
is actually additive, our approach can be modified to estimate the
parametric and nonparametric components easily. To estimate $\beta$ and
the two functions $\mu_t(t)$ and $\mu_z(\beta_0^TZ(t))$, we can perform
a two-step procedure. First, apply a one-dimensional scatter plot
smoother to $\{(Y_{ij},T_{ij})|i=1,\ldots,n;j=1,\ldots,N_i\}$ to
estimate $\mu_t(t)$. Then, apply modified rMAVE (\ref{rM2}) to the
residuals to estimate $\beta_0$. Specifically, $\beta_0$ can be
estimated by solving the minimization problem
\[
\min_{\beta,\mathbf{a},\mathbf{b},\mathbf{d}}\Biggl(\sum_{j=1}^n\sum_{\ell=1}^{N_j}
\sum_{i=1}^n\sum_{k=1}^{N_i}
[Y^c_{ik}-d_{j\ell}\beta^T(Z_{ik}-Z_{j\ell})]^2w_{ikj\ell}
\Biggr),
\]
where $Y^c_{ik} = Y_{ik}-\hat{\mu}_t(T_{ik})$ for $1\leq k \leq N_i$
and $1\leq i \leq n$.

Model (1.3) extends the popular single-index model from independent
univariate to longitudinal response data. Our extension allows both
time-independent and longitudinal covariates, but restricts the effect
of the covariates to be time invariant. Such a time-invariant approach
is in line with the philosophy in linear mixed-effects model, where the
covariate effect is time invariant. In this regard, our approach could
be viewed as an extension of a parametric linear mixed-effect model to
a more flexible semiparametric mixed-effects model.
While such an extension may still be considered restrictive, as
compared to an approach that adopts a time-dependent direction $\beta
(t)$ to model the covariate effects, the time-invariant direction $\beta
$ in our model has a nice interpretation as the averaging covariate
effect over time. Thus, $\beta$ as an average of $\beta(t)$ serves as a
simple summary measure for the (possibly more complicated)
time-dependent covariate effects. Moreover, a time-dependent approach
would require a lot more data to correctly estimate the direction $\beta
(t)$. When circumstances allow, one way to extend our approach to
time-dependent direction $\beta(t)$ is to adopt a two-stage procedure:
at the first stage, one bins the data in the direction of time and then
applies rMAVE to data that are observed in a bin that contains $t$ to
obtain an initial estimate of $\beta(t)$; these are smoothed at the
second stage to improve over the initial estimates. This is a subject
for future investigation and is beyond the scope of this paper.

Thus far, we have focused on estimation of the unknown components in
the functional single-index model. A functional principal component
analysis (FPCA) could be added after the covariate-adjusted mean
function has been estimated. The mean-adjusted FPCA (mFPCA)
proposed in \citet{JianW10} can be used to reconstruct the random
trajectories. More specifically, we can first show that the
asymptotic properties of the covariance estimator in Theorem 3.5 of
\citet{JianW10} hold by exploiting the $\sqrt{n}$-consistency of
$\hat{\beta}$. Then, the eigenvalues and eigenfunctions
corresponding to the estimated covariance can be estimated by
solving the eigenequation, and the PACE method proposed in \citet
{YaoMW051} can be used to
estimate the principal component scores and to select the number of
components.

\begin{appendix}
%s7 ###
\section{Assumptions}\label{appA}

The following type of continuity, as defined in \citet{Yao07},
will be needed: a real function $f(x,y)\dvtx \mathbb{R}^{n+m}
\rightarrow\mathbb{R}$ is continuous on $A\subseteq\mathbb{R}^n$,
uniformly in $y\in\mathbb{R}^m$, if, given any $x\in A$ and
$\delta<0$, there exists a neighborhood of $x$ not depending
on~$y$, say $U(x)$, such that $|f(x',y)-f(x,y)|<\delta$ for all
$x'\in U(x)$ and $y\in\mathbb{R}^m$. Our proofs cover both
time-independent and time-dependent covariates with slightly
different assumptions and arguments. For both cases, we assume the
observation times $T_{ij}$ are i.i.d. with probability density function
$f(t)$. For
a random design with time-invariant covariates, we assume that
$(T_{ij},\tilde{Z}_i,Y_{ij})$ have the same distribution as
$(T,\tilde{Z},Y)$ with joint probability density function $f_3(t,\tilde
{z},y)$, but
dependency is allowed among observations from the same subject.
The joint probability density functions of $(T,\tilde{Z})$ and
$(T_1,T_2,\tilde{Z},Y_1,Y_2)$ are denoted as $f_2(t,\tilde{z})$
and $f_5(t_1,t_2,\tilde{z},y_1,y_2)$, respectively. If the
covariate is longitudinal, then we assume that
$(T_{ij},\tilde{Z}_{ij},Y_{ij})$ have the same distribution as
$(T,\tilde{Z},Y)$ with joint probability density function $f_3(t,\tilde
{z},y)$ and the
joint probability density function of $(T_1,T_2,\tilde{Z}_1,\tilde
{Z}_2,Y_1,Y_2)$ is
$f_6(t_1,t_2,\tilde{z}_1,\tilde{z}_2,y_1,y_2)$.

Below, we describe the various assumptions that appear in the
theorems.

\begin{enumerate}[A.6$'$]
\item[A.1] $h_t \asymp h_z \asymp h$, $h \rightarrow0$,
$n\mathbb{E}(N)h^2\rightarrow\infty$, $\mathbb{E}(N)h^2\rightarrow0$ and
$n\mathbb{E}(N)h^6< \infty$.

\item[A.1$'$] $h_t \asymp h_z \asymp h$, $h \rightarrow0$,
$n\mathbb{E}(N)h^2\rightarrow\infty$, $\mathbb{E}(N)h\rightarrow0$ and
$n\mathbb{E}(N)h^6< \infty$.

\item[A.2] The number of observations $N_i(n)$ for the $i$th subject is
a random variable with $N_i(n)\sim N(n)$, where $N(n)$ is a
positive integer-valued random variable such that
$\lim\sup_{n\rightarrow\infty}\mathbb{E}N(n)^2/[\mathbb{E}N(n)]^2
< \infty$ and $N_i(n)$, $i=1,\ldots, n$, are i.i.d.

\item[A.3] The conditional mean
$\mu(t,\beta^Tz)=\mathbb{E}(Y|T=t,\beta^TZ=\beta^Tz)$ and its
derivatives up to second order are continuous on $\{(t,\tilde{z})\}$
and its derivatives up to the third order are bounded for all
$\beta\dvtx|\beta-\beta_0|<\delta$, for a $\delta>0$.

\item[A.4] The joint probability density function $f_2(t,\tilde{z})$ and its
derivatives up to
third order are bounded, up to second order they are continuous on
$\{(t,\tilde{z})\}$ and $f_2(t,\tilde{z})>0$ is bounded away from
zero for all $\beta\dvtx|\beta-\beta_0|<\delta$, for a $\delta>0$.

\item[A.5] The joint probability density function $f_3(t,\tilde{z},y)$ and
its derivatives up to
second order exist and are continuous on $\{(t,\tilde{z})\}$,
uniformly in $y\in\mathbb{R}$ for all
$\beta\dvtx|\beta-\beta_0|<\delta$, for a $\delta>0$.

\item[A.6] $f_6(t_1,t_2,\tilde{z}_1,\tilde{z}_2,y_1,y_2)$ is continuous
on $\{(t_1,t_2,\tilde{z}_1,\tilde{z}_2)\}$, uniformly in
$(y_1,\break y_2)\in\mathbb{R}^2$ for all $\beta\dvtx|\beta-\beta_0|<\delta$,
for a $\delta>0$.

\item[A.6$'$] $f_5(t_1,t_2,\tilde{z},y_1,y_2)$ is continuous on
$\{(t_1,t_2,\tilde{z})\}$, uniformly in $(y_1,y_2)\in\mathbb{R}^2$
for all $\beta\dvtx|\beta-\beta_0|<\delta$, for a $\delta>0$.
\end{enumerate}

The bandwidth assumption A.1 and the assumption on the measurement
schedule A.2 are required to show that the usual local properties of
the kernel estimators hold for longitudinal or functional data in
the presence of within-subject correlation. Assumptions A.3--A.6 are
regularity conditions for joint probability density functions and the
mean function. These
regularity conditions, along with the bandwidth assumption A.1, are
needed for the consistency results. A.$1'$ and A.$6'$ are the
assumptions when the covariate variables are time invariant. Adopting
assumption~A.4 is common practice in the theory of regression
estimation to study estimators on sets bounded away from the
troublesome regions [e.g., \citet{Hall89}, \citet{HardHI93}
and \citet{XiaTLZ02}].

%s8 ###
\section{\texorpdfstring{Proofs of (3.3) and
(3.4)}{Proofs of (3.3) and (3.4)}}
\label{appB}

\vspace*{-8pt}

\begin{pf*}{Proof of (\ref{Deq1})}
Let $(\beta,B)$ be an orthogonal matrix and, by Lemma \ref{lem2},
we obtain
%
%e8.1 ###
\begin{eqnarray} \label{Deq2}\quad
&& \{n\mathbb{E}N\hat{f_2}(t,\tilde{z})\}^{-1}
\sum_{i=1}^n\sum_{j=1}^{N_i}K_h(T_{ij}-t,\tilde{Z}_{ij}-\tilde
{z})(Z_{ij}-z)(Z_{ij}-z)^T
\nonumber\\
&&\qquad = \{n\mathbb{E}N\hat{f_2}(t,\tilde{z})\}^{-1}
\sum_{i=1}^n\sum_{j=1}^{N_i}K_h(T_{ij}-t,\tilde{Z}_{ij}-\tilde{z})(\beta,B)
\pmatrix{ \beta^T \cr B^T }\nonumber\\[-8pt]\\[-8pt]
&&\qquad\quad\hspace*{105.9pt}{}\times (Z_{ij}-z)(Z_{ij}-z)^T
(\beta,B)
\pmatrix{\beta^T \cr B^T} \nonumber\\
&&\qquad=  (\beta,B)
\pmatrix{h_z^2 & F_\beta(t,z)Bh_z^2 \cr
B^TF_\beta^T (t,z)h_z^2 & B^TG(z)B + O_p(h^2)}
\pmatrix{\beta^T \cr B^T} + O_p(h^3),\nonumber
\end{eqnarray}
where $G(z)=\mathbb{E}\{(Z_{ij}-z)(Z_{ij}-z)^T\}$ and $F_\beta
(t,z) = \frac{\partial}{\partial
\tilde{z}}(f_2\nu_\beta^T(t,z))/f_2$.

Next, consider
\begin{eqnarray*}
D_n^\beta & = &
\{n\mathbb{E}N\}^{-1}\sum_{i=1}^n\sum_{k=1}^{N_i}d_\beta^2(T_{ik},Z_{ik})
\{n\mathbb{E}N\hat{f_2}(T_{ik},\tilde{Z_{ik}})\}^{-1}
\\
&&\hspace*{72.3pt}{} \times
\sum_{j=1}^n\sum_{\ell=1}^{N_j}K_h(T_{j\ell}-T_{ik},\tilde{Z}_{j\ell
}-\tilde{Z}_{ik})\\
&&\hspace*{115.4pt}{}\times(Z_{j\ell}-Z_{ik})(Z_{j\ell}-Z_{ik})^T
\\
& = &
\{n\mathbb{E}N\}^{-1}\sum_{i=1}^n\sum_{k=1}^{N_i}d_\beta
^2(T_{ik},Z_{ik})(\beta,B)
\\
&&\hspace*{71.9pt}{} \times
\pmatrix{h_t^2 & F_\beta(T_{ik},Z_{ik})Bh_z^2 \cr
B^TF_\beta^T (T_{ik},Z_{ik})h_z^2 & B^TG(Z_{ik})B + O_p(h^2)}
%&&\hspace*{71.9pt}{} \times
\pmatrix{\beta^T \cr B^T}\\
&&{} + O_p(h^3) \\
& = & (\beta,B)
\pmatrix{\mathbb{E}\biggl\{\biggl(\dfrac{\partial\mu}{\partial\tilde{z}^0}\biggr)^2\biggr\}
h_t^2 & \tilde{F}Bh_z^2
\vspace*{3pt}\cr
B^T\tilde{F}^Th_z^2 & 2B^T\tilde{G}B + O_p(h^2+\delta_\beta)}
\pmatrix{\beta^T \cr B^T}\\
&&{} +
O_p(h^3+\delta_\beta h^2 ),
\end{eqnarray*}
where $\tilde{F}=\mathbb{E}\{(\frac{\partial\mu}{\partial
\tilde{z}^0})^2F_\beta(T,Z)\}$,
$\tilde{G}=\frac{1}{2}\mathbb{E}\{(\frac{\partial\mu}{\partial
\tilde{z}^0})^2G(Z)\}$, the second equality follows from
(\ref{Deq2}) and the last equality follows from Lemma \ref{lem1}.

Using the formula for matrix inverse in block form and letting
$\tau=\mathbb{E}\{(\frac{\partial\mu}{\partial
\tilde{z}^0})^2\}h_t^2$ and $\tilde{G}^*=(B^T\tilde{G}B)^{-1}$, we
obtain
\begin{eqnarray*}
\{D_n^\beta\}^{-1} & = & (\beta,B)
\pmatrix{\dfrac{1}{\tau} & \dfrac{-1}{2\tau}\tilde{F}B\tilde{G}^*h_z^2
\vspace*{2pt}\cr
\dfrac{-1}{2\tau}\tilde{G}^*B^T\tilde{F}^Th_z^2 &
\dfrac{1}{2}\tilde{G}^*}
\pmatrix{\beta^T \cr B^T} +
O_p(h+\delta_\beta) \\
& = & \frac{\beta\beta^T}{\tau}
-\frac{h_z^2}{2\tau}(B\tilde{G}^*B^T\tilde{F}^T\beta^T+\beta\tilde
{F}B\tilde{G}^*B^T)
+ \frac{1}{2}B\tilde{G}^*B^T+ O_p(h+\delta_\beta) \\
& = & \frac{\beta_0\beta_0^T}{\tau} -\frac{h_z^2}{2\tau}(\tilde{G}^+
\tilde{F}^T\beta_0^T+\beta_0\tilde{F} \tilde{G}^+) +
\frac{1}{2}\tilde{G}^+ + O_p(h+\delta_\beta),
\end{eqnarray*}
where $\tilde{G}^+ = B_0(B_0^T\tilde{G}B_0)^{-1}B_0^T$ is the
Moore--Penrose inverse of $\tilde{G}$.
\end{pf*}
\begin{pf*}{Proof of (\ref{Upeq1})}
To prove
\[
\Upsilon= \tilde{G}(\beta-\beta_0)-(n\mathbb{E}N
)^{-1} \sum_{i=1}^n\sum_{k=1}^{N_i} \biggl\{
\nu_{\beta_0}(T_{ik},Z_{ik})\,\frac{\partial\mu}{\partial
\tilde{z}^0} \biggr\}\epsilon_{ik} +o_p(n^{-1/2}),
\]
we work on $\{Y_{ik}-a_\beta(T_{j\ell},Z_{j\ell}) -
b_\beta(T_{j\ell},Z_{j\ell})(T_{ik}-T_{j\ell}) -
d_\beta(T_{j\ell},Z_{j\ell})\beta_0^T(Z_{ik}-Z_{j\ell})\}$ first.
By Lemma \ref{lem3}, we have
%
%e8.2 ###
\begin{eqnarray} \label{Upeq2}\qquad
&&
Y_{ij}  -a_\beta(t,z) - b_\beta(t,z)(T_{ij}-t) -
d_\beta(t,z)\beta_0^T(Z_{ij}-z) \nonumber\\
&&\qquad = -\frac{\partial\mu}{\partial
\tilde{z}}\,\nu_\beta(t,z)\delta_\beta+ \epsilon_{ij}-
\tilde{\varepsilon}_{n,1} \nonumber\\
&&\qquad\quad{}-
\frac{\tilde{\varepsilon}_{n,2}}{h_t}(T_{ij}-t) -
\frac{\tilde{\varepsilon}_{n,3}}{h_z}\beta_0^T(Z_{ij}-z) \nonumber\\
&&\qquad\quad{} + \frac{1}{2}\frac{\partial^2 \mu}{\partial t^2}\bigl(
(T_{ij}-t)^2 -h_t^2 \bigr)+ \frac{\partial^2 \mu}{\partial t\,
\partial\tilde{z}^0}(T_{ij}-t)(\tilde{Z_{ij}}-\tilde{z})
\\
&&\qquad\quad{} +
\frac{1}{2}\frac{\partial^2 \mu}{\partial
(\tilde{z}^0)^2}\bigl((\tilde{Z_{ij}}-\tilde{z})^2 - h_z^2\bigr)\nonumber\\
&&\qquad\quad{} + o_p( h^2 + |\delta_\beta|+|\delta_\beta|h+
|\delta_\beta|h^2 + |\delta_\beta|^2 ) \nonumber\\
&&\qquad\quad{} + O_p(\Delta),\nonumber
\end{eqnarray}
where
$\Delta=\sum_{\alpha_1+\alpha_2=3}|T_{ij}-t|^{\alpha_1}|\tilde
{Z}_{ij}^0-\tilde{z}^0|^{\alpha_2}+|T_{ij}-t||(Z_{ij}-z)^T\delta_\beta|+
|Z_{ij}^T-z^T|^2(\delta_\beta+\delta_\beta^2)$. We will now
calculate the weighted sum of each term in
(\ref{Upeq0}). This leads to the following results.

(i)
Let
\begin{eqnarray*}
\Sigma_{0n} & = & (n\mathbb{E}N )^{-1}
\sum_{j=1}^n\sum_{\ell=1}^{N_j} \frac{\partial\mu/
\partial\tilde{z}}{n\mathbb{E}N}\,\frac{d_{\beta}(T_{j\ell},Z_{j\ell
})}{\hat{f}_2}\\
&&\hspace*{74.4pt}{} \times\sum_{i=1}^{n}\sum_{k=1}^{N_i}
K_h(T_{ik}-T_{j\ell},\tilde{Z}_{ik}-\tilde{Z}_{j\ell})\\
&&\hspace*{116.5pt}{}\times (Z_{ik}-Z_{j\ell
})\nu_\beta^T(T_{j\ell},Z_{j\ell})
\\
& = & (n\mathbb{E}N )^{-1}
\sum_{j=1}^n\sum_{\ell=1}^{N_j} \frac{\partial\mu/ \partial
\tilde{z}\,d_{\beta}(T_{j\ell},Z_{j\ell})}{\hat{f}_2}\\
&&\hspace*{74.4pt}{}\times\bigl(f_2\nu_\beta
(T_{j\ell},Z_{j\ell})+ O_p(h^2) \bigr)\nu_\beta^T(T_{j\ell},Z_{j\ell})
\\
& = & (n\mathbb{E}N )^{-1}
\sum_{j=1}^n\sum_{\ell=1}^{N_j} \frac{\partial\mu/ \partial
\tilde{z}}{\hat{f}_2}f_2\nu_\beta(T_{j\ell},Z_{j\ell})\\
&&\hspace*{74.2pt}{}\times
\nu_\beta
^T(T_{j\ell},Z_{j\ell})\biggl(
\frac{\partial\mu}{\partial\tilde{z}}+O_p(h+\delta_\beta)\biggr)
+ o_p(n^{-1/2}) \\
& = & \tilde{G}+o_p(n^{-1/2}).
\end{eqnarray*}

(ii)
If we let
\begin{eqnarray*}
N_n & = & (n\mathbb{E}N )^{-1}
\sum_{j=1}^n\sum_{\ell=1}^{N_j}
\frac{d_{\beta}(T_{j\ell},Z_{j\ell})}{n\mathbb{E}N \hat{f}_2}
\sum_{i=1}^{n}\sum_{k=1}^{N_i}
K_h(T_{ik}-T_{j\ell},\tilde{Z}_{ik}-\tilde{Z}_{j\ell})\\
&&\hspace*{166.2pt}{}\times(Z_{ik}-Z_{j\ell})
\epsilon_{ik},
\end{eqnarray*}
then
%
%e8.3 ###
\begin{eqnarray} \label{Nneq1}
N_n & = & (n\mathbb{E}N )^{-1}
\sum_{i=1}^n\sum_{k=1}^{N_i} \frac{\epsilon_{ik}}{n\mathbb{E}N}
\sum_{j=1}^{n}\sum_{\ell=1}^{N_j}
K_h(T_{ik}-T_{j\ell},\tilde{Z}_{ik}-\tilde{Z}_{j\ell})\nonumber\\
&&\hspace*{132.8pt}{} \times\frac{d_{\beta
}(T_{j\ell},Z_{j\ell})(Z_{ik}-Z_{j\ell})}{\hat{f}_2}\\
& = & (n\mathbb{E}N )^{-1} \sum_{i=1}^n\sum_{k=1}^{N_i}
N_{nj}\epsilon_{ik},\nonumber
\end{eqnarray}
where
%
%e8.4 ###
\begin{eqnarray}\label{Nnjeq}\qquad
N_{nj} & = & \frac{1}{n\mathbb{E}N }
\sum_{j=1}^{n}\sum_{\ell=1}^{N_j}
K_h(T_{ik}-T_{j\ell},\tilde{Z}_{ik}-\tilde{Z}_{j\ell})\frac{d_{\beta
}(T_{j\ell},Z_{j\ell})(Z_{ik}-Z_{j\ell})}{\hat{f}_2}
\nonumber\\
& = & \frac{1}{n\mathbb{E}N} \sum_{j=1}^{n}\sum_{\ell=1}^{N_j}
K_h(T_{ik}-T_{j\ell},\tilde{Z}_{ik}-\tilde{Z}_{j\ell})(Z_{ik}-Z_{j\ell
})\nonumber\\
&&\hspace*{56.6pt}{}\times\frac{({\partial
\mu}/{\partial\tilde{z}}) +
O_p( h +|\delta_\beta| )}{f_2} \\
& = & \bigl(Z_i-\mathbb{E}(Z|T=T_{ik},Z^T\beta=Z_{ik}^T\beta)\bigr)\,
\frac{\partial \mu}{\partial\tilde{z}^0}+ O_p( h +|\delta_\beta| )
\nonumber\\
& = & -\nu_\beta(T_{ik},Z_{ik})\,\frac{\partial
\mu}{\partial\tilde{z}^0}+ O_p( h +|\delta_\beta| ).\nonumber
\end{eqnarray}
%
%d_\beta(t,z) = \frac{\partial\mu}{\partial\tilde{z}} +
% O_p( h +|\delta_\beta| )
Plugging (\ref{Nnjeq}) into (\ref{Nneq1}), we will get
\begin{eqnarray*}
N_n & = & (n\mathbb{E}N )^{-1}
\sum_{i=1}^n\sum_{k=1}^{N_i} \biggl(
-\nu_\beta(T_{ik},Z_{ik})\,\frac{\partial\mu}{\partial
\tilde{z}^0}+
O_p( h +|\delta_\beta| ) \biggr)\epsilon_{ik} \\
& = & (n\mathbb{E}N )^{-1} \sum_{i=1}^n\sum_{k=1}^{N_i}
\biggl( -\nu_{\beta_0}(T_{ik},Z_{ik})\,\frac{\partial\mu}{\partial
\tilde{z}^0} \biggr)\epsilon_{ik} + o_p(n^{-1/2}).
\end{eqnarray*}

(iii) From the definitions of
$\tilde{\varepsilon}_{n,1}$, $\tilde{\varepsilon}_{n,2}$ and
$\tilde{\varepsilon}_{n,3}$, we obtain
\[
\tilde{\varepsilon}_{n,1} = O_p(h^2), \qquad \tilde{\varepsilon}_{n,2} =
O_p(h^2)\quad\mbox{and}\quad\tilde{\varepsilon}_{n,3} = O_p(h^2).
\]
Let $R_n= (n\mathbb{E}N )^{-1}
\sum_{j=1}^n\sum_{\ell=1}^{N_j}
\frac{1}{n\mathbb{E}N}\sum_{i=1}^{n}\sum_{k=1}^{N_i}
K_h(T_{ik}-T_{j\ell},\tilde{Z}_{ik}-\tilde{Z}_{j\ell})
\times(Z_{ik}-Z_{j\ell}) ( \tilde{\varepsilon}_{n,1} +
\frac{\tilde{\varepsilon}_{n,2}}{h_t}(T_{ik}-T_{j\ell}) +
\frac{\tilde{\varepsilon}_{n,3}}{h_z}\beta_0^T(Z_{ik}-Z_{j\ell})
)$ and thus $R_n = o_p(n^{-1/2})$.

(iv)
\begin{eqnarray*}
&& (n\mathbb{E}N )^{-1} \sum_{j=1}^n\sum_{\ell=1}^{N_j}
\frac{\partial^2 \mu/ \partial t^2}{n\mathbb{E}N}\,\frac{d_{\beta
}(T_{j\ell},Z_{j\ell})}{\hat{f}_2}\\
&&\quad{} \times\sum_{i=1}^{n}\sum_{k=1}^{N_i}
K_h(T_{ik}-T_{j\ell},\tilde{Z}_{ik}-\tilde{Z}_{j\ell})(Z_{ik}-Z_{j\ell})
\{
(T_{ik}-T_{j\ell})^2 -h_t^2 \} \\
&&\qquad= (n\mathbb{E}N )^{-1}
\sum_{j=1}^n\sum_{\ell=1}^{N_j} h_t^2 \,\frac{\partial^2
\mu}{\partial t^2}\,
\frac{d_{\beta}(T_{j\ell},Z_{j\ell})}{\hat{f}_2}\times O_p(h^2) =
o_p(n^{-1/2}).
\end{eqnarray*}

\mbox{}\phantom{i}(v)
\begin{eqnarray*}
&& (n\mathbb{E}N )^{-1} \sum_{j=1}^n\sum_{\ell=1}^{N_j}
\frac{\partial^2 \mu/ \partial t\, \partial
\tilde{z}^0}{n\mathbb{E}N}\,\frac{d_{\beta}(T_{j\ell},Z_{j\ell})}{\hat
{f}_2}\\
&&\quad{} \times\sum_{i=1}^{n}\sum_{k=1}^{N_i}
K_h(T_{ik}-T_{j\ell},\tilde{Z}_{ik}-\tilde{Z}_{j\ell})(Z_{ik}-Z_{j\ell})
(T_{ik}-T_{j\ell})\beta^T(Z_{ik}-Z_{j\ell}) \\
&&\qquad= o_p(n^{-1/2}).
\end{eqnarray*}

(vi)
\begin{eqnarray*}
&& (n\mathbb{E}N )^{-1} \sum_{j=1}^n\sum_{\ell=1}^{N_j}
\frac{\partial^2 \mu/ \partial\tilde{z}^2}{n\mathbb{E}N}\,\frac{d_{\beta
}(T_{j\ell},Z_{j\ell})}{\hat{f}_2}\\
&&\quad{}\times\sum_{i=1}^{n}\sum_{k=1}^{N_i}
K_h(T_{ik}-T_{j\ell},\tilde{Z}_{ik}-\tilde{Z}_{j\ell})(Z_{ik}-Z_{j\ell})
[\{\beta^T(Z_{ik}-Z_{j\ell})\}^2 -h_z^2 ] \\
&&\qquad= o_p(n^{-1/2}).
\end{eqnarray*}

From (i)--(vi), the weighted sum of (\ref{Upeq2}) becomes
\begin{eqnarray*}
\Upsilon &=& \tilde{G}(\beta-\beta_0)+N_n+R_n+o_p(n^{-1/2}) \\
& = & \tilde{G}(\beta-\beta_0)-(n\mathbb{E}N )^{-1}
\sum_{i=1}^n\sum_{k=1}^{N_i} \biggl(
\nu_{\beta_0}(T_{ik},Z_{ik})\,\frac{\partial\mu}{\partial
\tilde{z}^0} \biggr)\epsilon_{ik} +o_p(n^{-1/2})
\end{eqnarray*}
and (\ref{Upeq1}) is thus proved.\vadjust{\goodbreak}
\end{pf*}

%s9 ###
\section{\texorpdfstring{Proof of Theorem \protect\lowercase{\ref{thm02}}}{Proof of Theorem 3.2}}
\label{appC}

In this Appendix, we consider $Y_{ij}$ as the $j$th
observation of $i$th subject made at a random time $T_{ij}$ with a
univariate longitudinal covariate $Z_{ij}$, where $i=1,\ldots,n$ and
$j=1,\ldots,N_i$. The following definitions are needed to derive the
asymptotic normalities of two-dimensional scatter plot smoothers.

A two-dimensional kernel function $K_2\dvtx\mathbb{R}^2\rightarrow
\mathbb{R}$ is of order $(\bolds{\nu},\kappa)$ if
%
%e9.1 ###
\begin{eqnarray} \label{eqK2}
&&\iint u^{k_1}v^{k_2}K_2(u,v)\,du\,dv \nonumber\\[-8pt]\\[-8pt]
&&\qquad=
\cases{0, &\quad $0\leq k_1+k_2<\kappa,k_1\neq\nu_1, k_2\neq\nu_2$, \cr
(-1)^{|\bolds{\nu}|}|\bolds{\nu}|!, &\quad $k_1=\nu_1, k_2=\nu_2$, \cr
\neq0, &\quad $k_1+k_2=\kappa$,}\nonumber
\end{eqnarray}
where $\bolds{\nu}$ is a multi-index $\bolds{\nu}=(\nu_1,\nu_2)$ and
$|\bolds{\nu}|=\nu_1+\nu_2$. Also, define the inverse Fourier
transform of $K_2(u,v)$ by
\[
\zeta_1(t,z)=\iint
\exp\bigl(-(iut+iwz)\bigr)K_2(u,w)\,du\,dw.
\]

Further, given an integer $Q\geq1$ and for $q=1,\ldots,Q$, let
$\psi_q\dvtx \mathbb{R}^3\rightarrow\mathbb{R}$ satisfy:

\begin{enumerate}[B.2]
\item[B.1] $\psi_q(t,z,y)$'s are continuous on $\{(t,z)\}$, uniformly in
$y\in\mathbb{R}$;%\label{B1}

\item[B.2] the functions $\frac{\partial^p}{\partial t^{p_1}\,\partial
z^{p_2}}\psi_q(t,z,y)$ exist for all arguments $(t,z,y)$ and are
continuous on $\{(t,z)\}$, uniformly in $y\in\mathbb{R}$, for
$p_1+p_2=p$ and $0\leq p_1,p_2 \leq p$.%\label{B2}
\end{enumerate}

The kernel-weighted averages for two-dimensional smoothers are
defined as
%
%e9.2 ###
\begin{equation}\label{gwa2}\hspace*{30pt}
\Psi_{q n}=\frac{1}{n\mathbb{E}Nh_t^{\nu_1+1}h_z^{\nu_2+1}}\sum_{i=1}^{n}\sum
_{j=1}^{N_i}\psi_q(T_{ij},Z_{ij},Y_{ij})
K_2\biggl(\frac{t-T_{ij}}{h_t},\frac{z-Z_{ij}}{h_z}\biggr),
\end{equation}
where $K_2$ is a kernel function of order $(\nu,\kappa)$ and $h_t$
and $h_z$ are bandwidths associated with $t$ and $z$, respectively.
The property of asymptotic normality of the local linear estimator
$\hat{\mu}(t,z)$ can be shown by using four specific $\psi_q$
functions. Let
\[
\alpha_q(t,z)=\frac{\partial^{|\bolds{\nu}|}}{\partial
t^{\nu_1}\,\partial z^{\nu_2}}\int\psi_q(t,z,y)f_3(t,z,y)\,dy
\]
and
\[
\sigma_{qr}(t,z)=\int
\psi_q(t,z,y)\psi_r(t,z,y)f_3(t,z,y)\,dy\|K_2\|^2,
\]
where $f_3(t,z,y)$ is the joint density of $(T,Z,Y)$,
$\|K_2\|^2=\int K_2^2$ and $1\leq q,r \leq Q$.
\begin{lem}\label{lem0}
Under assumptions \textup{A.2--A.6}, \textup{B.1--B.2}, $h_t\asymp h_z\asymp h$,
$h\rightarrow0$, $n\mathbb{E}(N)h^{|\nu|+2}\rightarrow\infty$,
$\mathbb{E}(N)h^2
\rightarrow0$ and $n\mathbb{E}(N)h^{2\bolds{\kappa}+2}< \infty$,
\[
\sqrt{n\mathbb{E}Nh_t^{2\nu_1+1}h_z^{2\nu_2+1}}
[(\Psi_{1n},\ldots,\Psi_{Qn})^T-(\mathbb{E}\Psi_{1n},\ldots,\mathbb{E}\Psi
_{Qn})^T]
\stackrel{\mathscr{D}}{\rightarrow}N(0,\Sigma).
\]
\end{lem}
\begin{pf}
This lemma can be shown by following similar procedures as used to
prove Lemma
C.1 in \citet{JianW10}. The only difference is in the
change-of-variable step of showing that $Q_2=o(1)$.
\end{pf}

The following two lemmas can be justified easily by following the
procedures in \citet{JianW10} and thus we omit the proof.
\begin{lem}\label{lem01}
Let $H\dvtx\mathbb{R}^Q \rightarrow\mathbb{R}$ be a function with
continuous first order derivatives,
$DH(v)=(\frac{\partial}{\partial
x_1}H(v),\ldots,\frac{\partial}{\partial x_Q}H(v))^T$ and
$\bar{N}=\frac{1}{n}\sum_{i=1}^{n}N_i$. Under assumptions \textup{A.2--A.6},
\textup{B.1--B.2}, $h_t\asymp h_z\asymp h$, $h\rightarrow0$,
$n\mathbb{E}(N)h^{|\nu|+2}\rightarrow\infty$, $\mathbb{E}(N)h^2
\rightarrow0$,
$\frac{h_z}{h_t}\rightarrow\rho_\mu$ and $
n\mathbb{E}(N)h_t^{2\bolds{\kappa}+2}\rightarrow\tau_\mu^2$ for some $0<
\rho_\mu, \tau_\mu< \infty$,
\begin{eqnarray*}
&&\sqrt{n\bar{N}h_t^{2\nu_1+1}h_z^{2\nu_2+1}}[H(\Psi_{1n},\ldots,\Psi
_{Qn})-H(\alpha_1,\ldots,\alpha_Q)]
\\
&&\qquad\stackrel{\mathscr{D}}{\rightarrow}
N(\beta_H,[DH(\alpha_1,\ldots,\alpha_Q)]^T \Sigma
[DH(\alpha_1,\ldots,\alpha_Q)]),
\end{eqnarray*}
where
$\Sigma=(\sigma_{qr})_{1\leq q, r\leq l}$ and
\begin{eqnarray*}\beta_H & = &
\sum_{k_1+k_2=\bolds{\kappa}}\frac{(-1)^{\bolds{\kappa}}}{k_1!k_2!}\biggl[
\int s_1^{k_1} s_2^{k_2} K_2(s_1,s_2)\,ds_1\,ds_2 \biggr] \\
&&\hspace*{32.5pt}{} \times\Biggl\{ \sum_{q=1}^Q
\frac{\partial H}{\partial\alpha_q}[(\alpha_1,\ldots,\alpha_Q)^T]\,
\frac{\partial^{k_1+k_2-\nu_1-\nu_2}}{\partial
t^{k_1-\alpha_q}\,
\partial z^{k_2-\nu_2}}\alpha_q(t,z)
\Biggr\}\tau_\mu\sqrt{\rho_\mu^{2k_2+1}}.
\end{eqnarray*}
\end{lem}
\begin{lem}\label{lem02}
Under the same assumptions as Lemma \ref{lem01}, together with the
assumption that the inverse
Fourier transform $\zeta_1(t,z)$ is absolutely integrable,
\[
{\sup_{t\in\mathcal{T};z\in\mathcal{Z}}}|\Psi_{q
n}-\alpha_q|=O_p\biggl(\frac{1}{\sqrt{n}h^{|\nu|+2}}\biggr)\qquad\mbox{where
}h\asymp h_t\asymp h_z.
\]
\end{lem}
\begin{pf*}{Proof of Theorem \ref{thm02}}
The theorem can be justified easily by employing Lemmas \ref{lem01} and
\ref{lem02}, and following the procedures used to prove Theorem \ref{thm02}
in \citet{JianW10}.
\end{pf*}

\vspace*{-12pt}

%s10 ###
\section{Auxiliary results}\label{appD}

\begin{lem}\label{lem1}
Suppose $\{\mathbf{T}_i,\mathbf{Z}_i,\mathbf{Y}_i\}$ are from an i.i.d.
sample, where $\mathbf{T}_i=(T_{i1},\ldots,T_{iN_i})$,
$\mathbf{Z}_i=(Z_{i1},\ldots,Z_{iN_i})$ and
$\mathbf{Y}_i=(Y_{i1},\ldots,Y_{iN_i})$. Let $\psi_s(t,z,y)$ be a
series of functions and assume that $\mathbb{E}\{\psi_s(T,Z,Y)\}$ and
$\operatorname{var}\{\psi_s(T,Z,Y)\}$ are both finite. Let
$\psi^T=(\psi_1,\ldots,\psi_p)$ and $\Psi^T = (\Psi_1,\ldots,\Psi_p)$,
where $ \Psi_s=
\frac{1}{n\mathbb{E}N}\sum_{i=1}^{n}\sum_{k=1}^{N_i}\psi_s(T_{ik},Z_{ik},Y_{ik})
$ for $s=1,\ldots,p$. Under assumptions \textup{A.1--A.6}, we obtain
%
%e10.1 ###
\begin{equation}\label{eq1}
\sqrt{n}\{\Psi-\mathbb{E}(\Psi)\}\stackrel{\mathscr
{D}}{\rightarrow}N_p(0,\Sigma),
\end{equation}
where
\begin{eqnarray*}
\Sigma &=& \frac{1}{\mathbb{E}(N)}\mathbb{E}\{\psi(T,Z,Y)\psi
^T(T,Z,Y)\}\\
&&{} + \frac{\mathbb{E}(N)-1}{\mathbb{E}(N)}\mathbb{E}\{\psi(T,Z,Y)\psi
^T(T',Z',Y')\}\\
&&{} - \mathbb{E}\{\psi(T,Z,Y)\}\mathbb{E}\{\psi^T(T,Z,Y)\}.
\end{eqnarray*}
Equation (\ref{eq1}) implies that
\[
\frac{1}{n\mathbb{E}N}\sum_{i=1}^{n}\sum_{k=1}^{N_i}\psi(T_{ik},Z_{ik},Y_{ik})
= \mathbb{E}\{\psi(T,Z,Y)\}+O_p(n^{-1/2}).
\]
\end{lem}
\begin{pf}
We can prove (\ref{eq1}) by showing that $\sqrt{n}\{a^T\Psi-
a^T\mathbb{E}(\Psi)\} \stackrel{\mathscr{D}}{\rightarrow}N_p(0,a^T\Sigma
a)$, where $a^T=(a_1,\ldots,a_p)$, by the central limit theorem.
\end{pf}
\begin{lem}\label{lem2}
Let
\[
\Phi_\beta(t,\tilde{z}) =
\frac{1}{n\mathbb{E}N}\sum_{i=1}^n\sum_{j=1}^{N_i}
K_h(T_{ij}-t,\tilde{Z}_{ij}-\tilde{z}) \biggl(
\frac{T_{ij}-t}{h_t}\biggr)^{\alpha_1}
\biggl(\frac{\tilde{Z}_{ij}-\tilde{z}}{h_z}
\biggr)^{\alpha_2}Y_{ij}.
\]
Suppose that $\mathbb{E}(Y|T=t,Z^T\beta=\tilde{z})=m(t,\tilde{z})$
and that assumptions \textup{A.1--A.6} hold. Then,
\begin{eqnarray*}
\Phi_\beta(t,\tilde{z}) & = &
m(t,\tilde{z})f_2(t,\tilde{z})\xi_{\alpha_1,\alpha_2} +
\frac{\partial}{\partial
t}\{m(t,\tilde{z})f_2(t,\tilde{z})\}\xi_{\alpha_1+1,\alpha_2}h_t
\\
&&{} + \frac{\partial}{\partial
\tilde{z}}\{m(t,\tilde{z})f_2(t,\tilde{z})\}\xi_{\alpha
_1,\alpha_2+1}h_z
+ O_p(h^2) \biggl(\mbox{or }O_p\biggl(\frac{1}{\sqrt{n\mathbb{E}N
h^2}}\biggr)\biggr),
\end{eqnarray*}
where $ \xi_{\alpha_1,\alpha_2}=\int
K(u,v)u^{\alpha_1}v^{\alpha_2}\,du\,dv$.
\end{lem}
\begin{pf}
From the definition of expectation and by the techniques of
change-of-variables and Taylor's expansion, we have
\begin{eqnarray*}
\mathbb{E}\{\Phi_\beta(t,\tilde{z})\} & = &
\int\frac{1}{h_t h_z}
K\biggl(\frac{s-t}{h_t},\frac{u-\tilde{z}}{h_z}\biggr) \biggl(
\frac{s-t}{h_t}\biggr)^{\alpha_1} \biggl(\frac{u-\tilde{z}}{h_z}
\biggr)^{\alpha_2}\\
&&\hspace*{9.3pt}{}\times yf_3(s,u,y)\,ds\,du\,dy \\
& = & \int
K(v_1,v_2)v_1^{\alpha_1}v_2^{\alpha_2}m(t+v_1h_t,\tilde
{z}+v_2h_z)\\
&&\hspace*{8.2pt}{}\times f_2(t+v_1h_t,\tilde{z}+v_2h_z)\,dv_1\,dv_2
\\
& = & m(t,\tilde{z})f_2(t,\tilde{z})\xi_{\alpha_1,\alpha_2} +
\frac{\partial}{\partial
t}\{m(t,\tilde{z})f_2(t,\tilde{z})\}\xi_{\alpha_1+1,\alpha_2}h_t
\\
&&{} + \frac{\partial}{\partial
\tilde{z}}\{m(t,\tilde{z})f_2(t,\tilde{z})\}\xi_{\alpha
_1,\alpha_2+1}h_z
+ O(h^2).
\end{eqnarray*}
The lemma now follows by Lemma \ref{lem0}.
\end{pf}

In the following lemma, we study the asymptotic expansions of the
weighted least-squares estimator,
$\hat{\theta}^T(t,z)=(\hat{a}_\beta(t,z),\hat{b}_\beta(t,z),\hat
{d}_\beta(t,z))$,
of the local linear smoother for mean function $\mu(t,\beta^Tz)$
when an initial single index $\beta$ is given. Thus,
%
%e10.2 ###
\begin{eqnarray} \label{thetaeq1}
\hat{\theta}(t,z) & = & \mathop{\arg\min}_\theta\sum_{i=1}^n\sum_{j=1}^{N_i}
K_h(T_{ij}-t,\tilde{Z}_{ij}-\tilde{z}) \nonumber\\
&&\hspace*{64.1pt}{} \times\bigl(Y_{ij}-a_\beta-
b_\beta(T_{ij}-t)-d_\beta(\tilde{Z}_{ij}-\tilde{z})\bigr)^2
\nonumber\\[-8pt]\\[-8pt]
& = &
\frac{1}{n\mathbb{E}N\Sigma_n^\beta(t,z)}\sum_{i=1}^n\sum_{j=1}^{N_i}
K_h(T_{ij}-t,\tilde{Z}_{ij}-\tilde{z}) \nonumber\\
&&\hspace*{94pt}{}\times\biggl(1,
\frac{T_{ij}-t}{h_t},\frac{\tilde{Z}_{ij}-\tilde{z}}{h_z}
\biggr)^TY_{ij},\nonumber
\end{eqnarray}
where
\begin{eqnarray*}
\Sigma_n^\beta(t,z) & = &
\frac{1}{n\mathbb{E}N}\sum_{i=1}^n\sum_{j=1}^{N_i}
K_h(T_{ij}-t,\tilde{Z}_{ij}-\tilde{z}) \\
&&\hspace*{57.2pt}{} \times\biggl(1,
\frac{T_{ij}-t}{h_t},\frac{\tilde{Z}_{ij}-\tilde{z}}{h_z}
\biggr)^T\biggl(1,
\frac{T_{ij}-t}{h_t},\frac{\tilde{Z}_{ij}-\tilde{z}}{h_z} \biggr).
\end{eqnarray*}
\begin{lem}\label{lem3}
Under assumptions \textup{A.1--A.6},
\begin{eqnarray*}
\hat{a}_\beta(t,z) & = & \mu(t,\tilde{z}^0) + \frac{\partial
\mu}{\partial\tilde{z}}\nu_\beta(t,z)\delta_\beta+
\frac{1}{2}\,\frac{\partial^2\mu}{\partial t^2}h_t^2 +
\frac{1}{2}\,\frac{\partial^2\mu}{\partial\tilde{z}^2}h_z^2 +
\tilde{\varepsilon}_{n,1} \\
&&{} + O_p( h^3 +
|\delta_\beta|h^2+ |\delta_\beta|h^3 +
|\delta_\beta|^2h ), \\
\hat{b}_\beta(t,z)h_t & = & \frac{\partial\mu}{\partial t}\,h_t +
\frac{\partial\mu}{\partial\tilde{z}}\, \frac{\partial
\nu_\beta(t,z)}{\partial t}\delta_\beta h_t +
\tilde{\varepsilon}_{n,2}\\
&&{} + O_p( h^3 +
|\delta_\beta|h^2+|\delta_\beta|h^3 +
|\delta_\beta|^2h ), \\
\hat{d}_\beta(t,z)h_z & = & \frac{\partial\mu}{\partial
\tilde{z}}\,h_z + \frac{\partial\mu}{\partial\tilde{z}}\,
\frac{\partial\nu_\beta(t,z)}{\partial\tilde{z}}\delta_\beta h_z
+ \tilde{\varepsilon}_{n,3}\\
&&{} + O_p( h^3 +
|\delta_\beta|h^2+|\delta_\beta|h^3 + |\delta_\beta|^2h ),
\end{eqnarray*}
where $\delta_\beta=\beta_0-\beta$, $\Sigma_n^\beta=
\Sigma_n^\beta(t,\tilde{z})$, $ \nu_\beta(t,z)=
\mathbb{E}(Z|T=t,Z^T\beta=z^T\beta)-z$, all the derivatives of
$\mu(t,\tilde{z})$ are evaluated at $(t,\tilde{z}^0)$ and
\[
\pmatrix{\tilde{\varepsilon}_{n,1} \cr
\tilde{\varepsilon}_{n,2} \cr
\tilde{\varepsilon}_{n,3}}
= (n\mathbb{E}Nf_2)^{-1}
\sum_{i=1}^n\sum_{j=1}^{N_i}K_h(T_{ij}-t,\tilde{Z}_{ij}-\tilde{z})
\pmatrix{1 \cr(T_{ij}-t)/h_t \cr(Z_{ij}-z)^T\beta/h_z}
\epsilon_{ij}.
\]
\end{lem}
\begin{pf}
By Lemma \ref{lem2}, we obtain
\begin{eqnarray*}
\Sigma_n^\beta(t,z) &=&
\pmatrix{f_2 & \dfrac{\partial f_2}{\partial t}\,h_t &
\dfrac{\partial f_2}{\partial\tilde{z}}\,h_z \vspace*{2pt}\cr
\dfrac{\partial f_2}{\partial t} \,h_t & f_2 & 0 \vspace*{2pt}\cr
\dfrac{\partial f_2}{\partial\tilde{z}}\,h_z & 0 & f_2}
+ O_p(h^2),\\
\det\{\Sigma_n^\beta(t,z)\}&=&(f_2)^3+O_p(h^2)
\end{eqnarray*}
and
\begin{eqnarray*}
\{ \Sigma_n^\beta(t,z)\}^{-1} &=&
\frac{1}{f_2(t,\tilde{z})} \left( I -
\frac{1}{f_2(t,\tilde{z})}\pmatrix{0 & -\dfrac{\partial
f_2}{\partial t}\,h_t & -\dfrac{\partial f_2}{\partial
\tilde{z}}\,h_z \vspace*{2pt}\cr
-\dfrac{\partial f_2}{\partial t}\,h_t &
0 & 0 \vspace*{2pt}\cr
-\dfrac{\partial f_2}{\partial\tilde{z}}\,h_z & 0 & 0}
\right)\\
&&{} + O_p(h^2).
\end{eqnarray*}

By Taylor's expansion at $ (t,\tilde{z}^0) $,
$Y_{ij}=\mu(T_{ij},Z_{ij}^T\beta_0)+\epsilon_{ij}$ can be
expressed as
\begin{eqnarray*}
Y_{ij} & = & \underbrace{\mu(t,\tilde{z}^0) + \frac{\partial
\mu}{\partial t}\,(T_{ij}-t) + \frac{\partial\mu}{\partial
\tilde{z}^0}\,(\tilde{Z}_{ij}-\tilde{z})}_{E_1}{} + {}
\underbrace{\frac{\partial\mu}{\partial
\tilde{z}^0}\,(Z_{ij}-z)^T\delta_\beta}_{E_2}\\
&&{} +
{}\underbrace{\frac{1}{2}\,\frac{\partial^2 \mu}{\partial
t^2}\,(T_{ij}-t)^2}_{E_3}
+ \underbrace{\frac{\partial^2 \mu}{\partial t\,
\partial\tilde{z}^0}\,(T_{ij}-t)(\tilde{Z}_{ij}-\tilde{z})}_{E_4}\\
&&{} +
{}\underbrace{ \frac{1}{2}\,\frac{\partial^2 \mu}{\partial
(\tilde{z}^0)^2}\,(\tilde{Z}_{ij}-\tilde{z})^2}_{E_5}{} + {}O_p(E_6) +
\epsilon_{ij},
\end{eqnarray*}
where $E_6 =
\sum_{\alpha_1+\alpha_2=3}|T_{ij}-t|^{\alpha_1}|\tilde{Z}_{ij}^0-\tilde
{z}^0|^{\alpha_2}+
|T_{ij}-t||(Z_{ij}-z)^T\delta_\beta|+
|Z_{ij}^T-z^T|^2(|\delta_\beta|+|\delta_\beta|^2)$ and the
derivatives of $\mu(t,\tilde{z})$ are evaluated at
$(t,\tilde{z}^0)$. Therefore, $\hat{\theta}(t,z)$ in
(\ref{thetaeq1}) is the sum of the weighted averages of
$E_1,\ldots,E_6$ and error $\epsilon_{ij}$. After evaluating the
weighted average of each term, which amounts to smoothing each $E_i$
and $\epsilon_{ij}$, the lemma follows by combining all of the
smoothing terms.
\end{pf}
\end{appendix}

\section*{Acknowledgments}
The authors would like to thank the Editor, Associate Editor and two
referees for insightful comments.

%suskaldyti doi

%
\printaddresses

\end{document}